\documentclass[11pt]{article}
\usepackage{pst-plot,pstricks,pst-node}
\usepackage[utf8]{inputenc}
\usepackage{xspace}
\usepackage{epsfig,graphicx}
\usepackage{amsmath,amssymb}
\usepackage{theorem}
\usepackage{setspace}
\usepackage{multirow}
\usepackage{fancyhdr} 
\usepackage[numbers,sort]{natbib} 
\usepackage{empheq}
\usepackage{cleveref}

\usepackage{tikz}
\usepackage{pgfplots}
\usepackage{cancel}
\usepackage{enumerate}
\usepackage{rotating}
\usepackage{pdflscape}
\usepackage{afterpage}
\usepackage{longtable}
\usepackage{url}

\onehalfspace

\usepackage[lined,algoruled,linesnumbered,titlenumbered]{algorithm2e} 
\allowdisplaybreaks[1]

\SetKwRepeat{Do}{do}{while}%

\setcitestyle{authoryear,open={(},close={)}}

\fancypagestyle{memo}{
	\pagestyle{fancy}
	\fancyhf{}%
	\fancyhead[RO]{\thepage}%
	\fancyhead[LO]{\rightmark}%
	
}%

\theoremstyle{plain}
\newtheorem{prop}{Proposition}
\newtheorem{lemma}{Lemma}

\theoremstyle{plain} {
	\theorembodyfont{\rmfamily}%

}

\newcommand{\ProofNoNL}{{\bf \noindent Proof.}\xspace}

\newcommand{\EndProofNoNL}{\hfill $\Box$ \par \bigskip}

\newcommand{\wloge}{without loss of generality\xspace}

\newcommand{\ox}{{\ensuremath{\overline X}}\xspace}
\newcommand{\og}{{\ensuremath{\overline \gamma}}\xspace}

\makeatletter
\renewcommand*\env@matrix[1][c]{\hskip -\arraycolsep
	\let\@ifnextchar\new@ifnextchar
	\array{*\c@MaxMatrixCols #1}}
\makeatother 

\DeclareMathOperator*{\argmax}{arg\,max}

\newcommand{\name}{\mbox{D-MCLP}\xspace} 
\newcommand{\namef}{\mbox{BL-D-MCLP}\xspace} 
\newcommand{\nameat}{\mbox{A-MCLP}\xspace} 

\usepackage[normalem]{ulem}

\newcommand{\msout}[1]{\text{\sout{\ensuremath{#1}}}}
\definecolor{armygreen2}{rgb}{0.99, 0.53, 0.93}
\newcommand{\AMR}[2]{{\color{armygreen2}#1 \ifmmode\msout{#2}\else\sout{#2}\fi}}

\definecolor{verde}{rgb}{0,0.4,0}
\definecolor{verdeclaro}{rgb}{0.01, 0.75, 0.24}
\title{Edge downgrades in the Maximal Covering Location Problem} 
\author{\smaller Marta Baldomero-Naranjo$^1$, 
J\"org Kalcsics$^2$, Antonio M. Rodr\'iguez-Ch\'ia$^1$ \\[1ex]
\smaller $^1$ Departamento de Estad\'istica e Investigaci\'on Operativa, Universidad de C\'adiz, Facultad de Ciencias, \\ \smaller    11510, Puerto Real (C\'adiz), Spain,  marta.baldomero@uca.es, antonio.rodriguezchia@uca.es\\
\smaller $^2$ School of Mathematics and Maxwell Institute for Mathematical Sciences, \\ \smaller University of Edinburgh, Edinburgh, EH9 3FD, United Kingdom, joerg.kalcsics@ed.ac.uk 
}
\date{}

\begin{document}

\maketitle

\begin{abstract}
We tackle the downgrading maximal covering location problem within a network. In this problem,   
two actors with conflicting objectives are involved:
a) The location planner aims to determine the location of facilities to maximise the covered demand while anticipating that an attacker will attempt to reduce coverage by increasing the length of some edges (downgrade); 
b) The attacker seeks to maximise the demand initially covered by the facilities but left uncovered after the downgrade. The attacker can increase the length of certain edges within a specified budget.

We introduce a mixed-integer linear bilevel program to formulate the problem, followed by a preprocessing phase and a matheuristic algorithm designed to address it. Additionally, computational results are presented to illustrate the potential and limitations of the proposed algorithm.
\end{abstract}


\section{Introduction}

In the area of facility location, the Maximal Covering Location Problem is a well-known classical problem that has attracted the attention of numerous researchers since its introduction by \cite{ChuRev74}. The majority of the studies that deal with this problem assume that the distances from the clients to the facilities remain unalterable, see e.g. \cite{CorFurLju19,BalKalRod20, BlaGazSal23, MarMarRodSal18,GarMar19, Mur16,BerDreKar10,BalMarRod24,FarAsg12} and references therein.

In this paper, we propose an extension that considers that an attacker can modify the length of some edges within a budget and up to a certain value to increase the distance from the clients to the facilities. Under this assumption, the Downgrading Maximal Covering Location Problem (\name) emerges as a notable challenge, involving the conflicting interests of both the facility location planner and potential attackers. At its core, the \name entails the strategic location of facilities within the nodes of a network to maximise coverage while anticipating and mitigating potential attacks aimed at subsequently reducing this coverage. This problem comprises an interaction between two distinct actors, each with conflicting objectives:

\begin{enumerate}
\item First, the location planner undertakes the task of strategically locating  facilities to ensure maximal coverage of demand nodes within the network after the attacker has acted. This entity operates in anticipation of adversarial actions, aware of the fact that the potential attacker seeks to optimally sabotage the effectiveness of the chosen facility locations by increasing edge lengths. 

\item Second, the attacker (they) seeks to minimise the demand still covered by the facilities after downgrading the network, i.e., maximise the demand initially covered by the facilities but left uncovered after the downgrade of the network. To achieve that, they can increase the length of some edges within a given budget and up to a certain value.
\end{enumerate} 

Therefore, first the location planner locates the facilities, and then the attacker downgrades the network. Note that this implies that the distances in the network are modified during the optimization problem, thus, the distance between two nodes after the downgrades will have to be computed within the model.

This model has notable applications across various domains, especially in the design of systems that account for potential disruptions. For example in disaster relief, supply depots need to be set up such that supplies can reach affected regions within a certain time. However, some roads may have become difficult to traverse, possibly making it impossible to reach some regions within reasonable time. Although the damage caused by natural disasters is not a deliberate action, planners may want to hedge against worst case scenarios when setting up supply depots. More direct applications arise in defense and security. In the former, military strategists will want to locate supply points behind the front line where additional troops can be stationed and material can be stored. In case of an attack, reinforcements need to be able to reach defensive frontline positions quickly and the attacker may decide to divert artillery or airborne units to slow down those reinforcements as best as possible. In wildlife conservation, the model can help determine the location of guard and patrol stations to protect wildlife, anticipating that poachers might downgrade certain paths (e.g., by blocking forest trails with trees, making them more difficult to traverse). 

This study falls within the domain of downgrading (upgrading) network problems, in which specific elements initially treated as fixed inputs in the classical version are now transformed into decision variables (e.g.~the length of the edges). These variables are then simultaneously determined by agents, whose
 task is to change the underlying network within certain limits such that the optimal objective value that can be obtained in
the modified network is as bad (good) as possible. In the following, we briefly review the related literature.

\subsection{Literature Review}
In the context of downgrading location problems, we are only aware of studies on node downgrades (the weight of the vertices can be modified subject to a prespecified budget): the 1-median problem \citep{Gassner07UpDown1median,Gas09}, the 1-center problem \citep{Gassner09UpDown1center}, and the Euclidean 1-median problem \citep{Plas16}. 
As far as we know, this is the first study analysing edge downgrades.

Regarding upgrading problems, several studies have been published recently in the context of either node upgrades or edge upgrades, for example the maximal covering \citep{BalKalMarRod22, BalKalRod24}, the $p$-median \citep{sepasian15uppmedian,AfrAliBar20,EspMar23}, the $p$-center \citep{Sepasian181center,AntLanSal23}, the hub-location \citep{BlaMar19,LanMunRodSal24}, 
the graphical TSP \citep{LanPla23}, the spanning tree \citep{ALVAREZMIRANDA201713}, the
maximal shortest path interdiction \citep{ZhaGuaPar21}, among others. 

The \name is also closely related to interdiction problems \citep{Smith2013,FisLjuMichSin19, SmiSon20, RamHaySinVid24,OhChu11}. While we propose increasing the length of one or more edges in the \name, interdiction problems usually consider the case of removing edges/arcs instead. The covering-interdiction problem \citep{FroRuz21}, in which the first player (the location planner) locates some facilities in a network to maximise the amount of demand covered by the facilities, and the second player, called the interdictor, disrupts the infrastructure by removing at most a given number of edges in the network to reduce the coverage as much as possible. Their paper is focused on analysing the problem's complexity.  
Observe that the interdiction action in this context is a particular case of downgrading, where the length of an interdicted edge is increased to a value greater than the covering radius. 
Similarly, given a subset of nodes chosen as facilities, the reachability-interdiction problem \citep{FroRuz21tree} studies the problem where an interdictor can remove a set of edges to prevent as many vertices as possible from reaching a facility. A node can no longer reach a facility after removing edges if no facility is located in the connected component it belongs to. 
Furthermore, \cite{CapSca11} propose a multilevel program for allocating protection resources among the components of a shortest-path network so as to maximize its robustness to external disruptions. Regarding the shortest path interdiction problem, \cite{Smith2013} considers the situation where the attacker seeks to maximize the length of a shortest $s$-$t$-path by either increasing the cost of some arcs or by deleting arcs altogether. 
A similar version of this shortest path interdiction problem appears as a component of our problem from the attacker's perspective. The most vital elements for the $p$-median and $p$-center problems \citep{bazTouVan13} are also slightly related since a fixed number of edges are removed to provide their worst objective values. Finally, it is worth mentioning the Maximal Covering Location Disruption Problem \citep{Lun24}, a zero-sum game where the leader seeks to disrupt the MCLP solution by making a subset of the possible
facility locations unavailable. 
Note that the design of robust networks against attacks has been studied from a game-theoretic approach for a wide range of problems, as for example in \cite{PerPue13}.

\subsection{Overview}

The main contribution of this paper are: 
\begin{itemize}
    \item The first model, (\name), to deal with continuous edge downgrades in the Maximal Covering Location Problem is introduced. A mixed-integer linear bilevel formulation is proposed which is the reformulation of an intuitive three-level formulation. 
    \item A strategy to preprocess the data to reduce the number of variables and constraints of the formulation is introduced. Moreover, strategies to improve the formulation of the attacker's problem given the location of the facilities are proposed. 
    
    \item A matheuristic algorithm to solve the problem is developed. For this algorithm, we present some variants (levels of intensity), which allow the user to decide the best trade-off between the computational time employed and the quality of the solution.
    \item  Computational experiments are carried out to 
    illustrate the increase of the covered demand after the attack using the \name rather than using the classical MCLP. Furthermore, we show the potential and limitations of the proposed matheuristic algorithm. The solutions provided by the matheuristic are compared with the ones obtained by the general purpose solver for bilevel programming problems developed in \cite{FisLjuMonSin17}. 
\end{itemize}

The contents of the paper are organized as follows. The formulation of \name is presented in Section~\ref{sec:ProblemDescription}.  The strategy to preprocess the data is introduced in Section~\ref{sec:Prep}.  The matheuristic is discussed in Section~\ref{sec:Heuristic}. In this algorithm, we need to solve the attacker's problem given the location of the facilities several times, for this reason, we also propose some strategies to improve its formulation in this section. In Section~\ref{sec:ComputationalResults}, we discuss the computational results. 
Our conclusion and future line of research are presented in Section~\ref{sec:Conclusion}. Finally, in Appendix~\ref{sec:App}, we include the proof that the attacker's problem is NP-hard for a fixed set of facility locations.

\section{Problem description and mathematical formulations} \label{sec:ProblemDescription}

Let $N=(V,E,\ell)$ be an undirected network with node set $V$ and edge set $E$. 
Every edge $e=[i,j]=[j,i]\in E$, $i,j \in V,$ has a positive length $\ell_e=\ell_{[i,j]}$. For $i,j\in V,$ $d(i,j)$ is the shortest path distance between $i$ and $j$ before downgrading. Moreover, we are given a fixed coverage radius $R > 0$. We say that a node $i\in V$ is \emph{covered} by a facility at node $j$ if it is within a distance of strictly less than the coverage radius $R$ from a facility, i.e., $d(i,j) < R.$ Observe that we are deviating slightly from the standard definition, where $i$ was covered by $j$ whenever $d(i,j)\leq R$; as will be discussed in Section~\ref{sec:mid_level_problem}, this change is necessary to ensure a closed feasible region of the attacker's problem. Additionally, we are given a non-negative \emph{weight} $w_i$ that specifies the demand at the node $i$, for $i\in V$.

The length $\ell_e$ of each edge $e\in E$ can be increased by an amount lower than or equal to $u_e$ (known in advance). Let $\gamma = (\gamma_e)_{e \in E}$ be a variable vector defining the length increases of each edge, then $0 \le \gamma_e \le u_e$, for $e\in E$. Furthermore, increasing the length of edge $e \in E$ by an amount $\gamma_e \in [0, u_e]$ comes at a positive cost of $c_e \gamma_e $, for $c_e \in \mathbb{R}^+$, and there is a budget constraint $B \in \mathbb{R}^+$ on the overall cost of increase. Moreover, let $d(i,j,\gamma)$ be the shortest path distance between $i$ and $j$ after the edge length increases $\gamma$ have been applied, i.e., the shortest path in the network $N(\gamma):=(V,E,\ell(\gamma))$ where $\ell_e(\gamma) = \ell_e + \gamma_e,$ for $e \in E$. 

Finally, for $p \in \mathbb{N},$ let $X \subseteq V$ denote a set of $p$ nodes, let $C(X) = \{i \in V \mid \exists j \in X: d(i,j) < R \}$ denote the set of all nodes covered by a facility in $X$ before the edge length increases,  let $\mathcal {C}(X,\gamma) = \{i \in C(X) \mid \exists j \in X: d(i,j,\gamma) < R \}$ denote the set of all nodes covered by a facility in $X$ after the edge length increases, and let $\mathcal {\overline {C}}(X,\gamma) = \{i \in C(X) \mid \forall j \in X : d(i,j,\gamma) \geq R \}$ denote the set of all nodes initially covered by a facility in $X$ but no longer covered after the edge length increases, i.e., $\mathcal {\overline {C}}(X,\gamma)=C(X) \setminus \mathcal {C}(X,\gamma)$.
Then, \name
can be formulated as:
\[ \max_{X \subseteq V, |X|=p} \left\{ \sum_{i \in \mathcal C(X,\gamma)}\, w_i\: \Big|\: \gamma\in \argmax \left\{ \sum_{i \in \mathcal {\overline {C}}(X,\gamma)}\, w_i\: \Big|\:\sum_{e \in E} c_e \gamma_e \le B, 0 \le \gamma_e \le u_e, e \in E\right\}  \right\}.
\]
Observe that, although the problem is presented for undirected networks, the formulation we propose requires the use of a directed network. Specifically, we define the arc set $A$ containing arcs $a_1=(i,j)$ and $a_2=(j,i)$ for each edge $[i,j]\in E$. We refer to this edge $[i,j]\in E$, as the underlying edge of an arc $a_1$ or $a_2$ and it is denoted by $e_{a_1}$ or $e_{a_2}$, i.e., $e_{a_i}=e_{a_2}=[i,j]$. Table~\ref{tb:notation} summarises the main notation presented before. 
\begin{table}[htb]
\caption{Notation used throughout the paper.}
\medskip

\begin{tabular}{lp{0.75\textwidth}}
\hline
$N = (V,E, \ell)$ & The network with the set of nodes $V,$  the set of edges $E,$ and edge length $\ell$. \\
$N(\gamma)=(V,E,\ell(\gamma))$ & The network with the set of nodes $V,$  the set of edges $E,$ and edge lengths $\ell_e(\gamma) = \ell_e + \gamma_e,$ for $e \in E$.  \\
$c_e$& The unit cost of increasing the length of edge $e,$ $e\in E.$ \\ 
$e_a$& The underlying edge of an arc $a$, e.g., for $(i,j)\in A, e_a=[i,j].$ \\
$\ell_e(\gamma)$ & The length of edge $e,$ $e\in E,$ after the increase $\gamma$, i.e., $\ell_e + \gamma_e$. \\
$u_e$& Maximum amount that the length of edge $e$ can be increased, $e\in E$.\\
$w_i$& The demand of node $i\in V$. \\ 
$p$& The number of facilities to locate. \\ 
$X$& A set of $p$ facility locations, $X \subset V$, $|X|=p$. \\ 
$R$ & The coverage radius.\\
$B$& The budget for downgrading edges.  \\
\hline 
\end{tabular}	

\label{tb:notation}
\end{table}
\medskip

In the next subsection, we present the different decision levels of the problem.

\subsection{Trilevel Problem}

\noindent
Our downgrading maximal covering location problem can be characterized as a trilevel max-max-min problem:

\begin{description}
\item[Upper level:] The defender or location planner wants to locate a given number $p$ of facilities to maximise the amount of covered demand after the attacker has acted. For modelling its decision, we introduce the following family of binary variables: 

\begin{tabular}{lp{0.78\textwidth}}

$x_j$ & 1, if a facility is located at node $j$, 0 otherwise, for $j\in V$.\\
\end{tabular}	

\item[Middle level:] For a given set $X\subset V$ of facility locations, the attacker wants to downgrade edges, i.e., increase their length subject to a given budget to maximise the amount of demand that was initially covered by the facilities and is no longer covered after the downgrade. 
 For modelling its decision, we introduce the following family of continuous variables: 
 
 \begin{tabular}{lp{0.78\textwidth}}
$\gamma_e=\gamma_{[i,j]}$&  The increase in the length of edge $e=[i,j]$, for $e\in E.$\\
\end{tabular}	

\item[Lower level:] After the attacker has downgraded some edges, the defender needs to check if a node is still covered by computing the length of a shortest path to the nearest facility in the downgraded network and comparing it to the coverage radius $R$. Strictly speaking, for this step it is not necessary to find the shortest path. It suffices to check whether there exists a path of length strictly smaller than $R$ or not.
\end{description}

Observe that, unless stated otherwise, a solution to an instance of the trilevel program only refers to the optimal solutions of the first two stages, as an optimal solution for the third stage can easily be computed from the ones for the first two and it is, hence, not the primary focus for the decision-makers.

\subsubsection{Lower level problem}
We start by modelling the lower level problem. Let the upper and middle level decisions be given and denoted as $\overline X$ and $\overline \gamma = (\overline \gamma_e)_{e \in E}$, respectively. To check if node $i\in C(\overline X)$ is (still) covered after downgrading, we compute the length of a shortest $i$-$t$-path in an auxiliary induced directed network $N^i(\overline X) = (V^i, A^i(\overline X), \ell(\overline \gamma))$, where $t$ is an auxiliary node, $V^i = V \cup \{t\}$ and $A^i(\overline X) = \{(k,l), (l,k) \mid [k,l] \in E \}  \cup \{(x,t) \mid x \in \overline X\}$. 
The length of an arc $a$ is identical to the length of its underlying edge $e_a$. 
The lengths of the remaining arcs are fixed to zero and the values of $\gamma$-variables for the corresponding underlying edges are also zero. 
Finally, for $k\in V^i$, let $\Gamma_k^{i,+} \subset A^i(\overline X)$ and $\Gamma_k^{i,-} \subset A^i(\overline X)$ be the set of outgoing and, respectively, incoming arcs for node $k \in V^i$ in $N^i(\overline X)$.

Before presenting the integer programming formulations, we introduce some additional decision variables:

\begin{table}[h]
\begin{tabular}{lp{0.78\textwidth}}
$f_a^i$& The amount of flow traversing arc $a$,
for $i\in V, a\in A^i(\overline X)$. \\[0.5ex]
\end{tabular}	
\end{table}

\noindent
Let $d_i(\overline X, \overline \gamma)$ be the length of a shortest path from node $i$ to the nearest facility in $\overline X$ after the downgrade $\overline \gamma$, for $i\in V$, i.e., the shortest path form $i$ to $t$ in $N^i(\overline X)$. Then, $d_i(\overline X, \overline \gamma)$ can be computed as
\begin{alignat}{4}
d_i(\overline X, \overline \gamma) \:=\: \min & \sum_{a\in A^i(\overline X)} f_a^i(\ell_{e_{a}} +\overline \gamma_{e_{a}}) & \\
\text{s.t.} &\sum_{a\in \Gamma_k^{i,+}} f_{a}^{i} - \sum_{a\in \Gamma_k^{i,-}} f_{a}^{i} & \:=\: \begin{cases} 1, & \text{ for } k = i, \\ -1, & \text{ for } k = t, \\ 0, & \text{ for } k \in V, \\ \end{cases} \qquad && \\
&f_{a}^{i}  \:\geq\: 0, && a \in A^i(\overline X),
\end{alignat} 
or, using the dual of the shortest path problem, as
\begin{alignat}{4}
d_i(\overline X, \overline \gamma) \:=\: \max\; \pi^i_{i} &- \pi^i_t \\
\text{s.t. } \pi^i_k &\le \pi^i_l + \ell_{e_a} + \overline \gamma_{e_a}, \qquad \qquad && a=(k,l) \in A^i(\overline X), \label{constr:Level3:Dual:NodePotential} \\
\qquad \qquad \pi^i_k & \in \mathbb{R,}
&& k \in V^i,
\end{alignat} 
where $\pi_k^i$ is the node potential of node $k$, for $k\in V^i$ \citep{FulHar77,Gol78}. 
\medskip

Since the location of the facilities is a decision variable in our problem, i.e., $\overline{X}$ is not known a priori, we first discuss how to formulate the shortest path problem without using the sets of arcs $A^i(\overline{X})$. 
To that end, we modify the network $N^i(\overline X)$ as follows. We define $N^i=(V^i,A^i, \ell(\overline \gamma))$, where $V^i$ remains unchanged and $A^i = \{(k,l), (l,k) \mid [l,k] \in E\} \cup \{(j,t) \mid j \in V\}$, i.e., each node is now connected to $t$. The length of an arc $a$ is again identical to the length of its underlying edge $e_a$, while the lengths of the remaining arcs are fixed to zero and the values of $\gamma$-variables for the corresponding underlying edges are also zero.
Then, $d_i(\overline{X}, \overline \gamma)$ can be computed as
\begin{alignat}{4}
d_i(\overline X, \overline \gamma) \:=\: \min &\sum_{a\in A^i} f_a^i(\ell_{e_{a}}+\overline \gamma_{e_{a}}) & \\
\text{s.t.} &\sum_{a\in \Gamma_k^{i,+}} f_{a}^{i} - \sum_{a\in \Gamma_k^{i,-}} f_{a}^{i}  \:=\: \begin{cases} 1, & \text{ for } k = i, \\ -1, & \text{ for } k = t, \\ 0, & \text{ for } k \in V, \\ \end{cases} \qquad && \label{constr:Level3:FlowBalance} \\
&f_{a}^{i} \:\leq\: \overline x_k, && a = (k,t) \in A^i,\label{constr:Level3:LocVariables} \\
&f_{a}^{i} \:\geq\: 0, && a \in A^i.
\end{alignat} 
Constraints \eqref{constr:Level3:LocVariables} ensure that we can not use arc $(k,t)$, $k \in V$, if there is no facility at node $k$. For the sake of understanding, the notation of the variables is the same as in the previous formulation. However, they are not exactly the same since their index domains differ, i.e, now the variables $f_a^i$ are defined for $a\in A^i$ while before they were defined in $a\in A^i(\overline X)$.  The dual is then given as
\begin{alignat}{4}
d_i(\overline X, \overline \gamma) \:=\: \max\; &\pi^i_{i} - \pi^i_t - \sum_{k \in V} \overline x_k \mu_k \label{constr:Level3:Dual2:ObjectiveFunction}\\
\text{s.t. } &\pi^i_k \le \pi^i_l + \ell_{e_a} + \overline \gamma_{e_a}, \qquad \qquad && a=(k,l) \in A^i, l \neq t, \label{constr:Level3:Dual2:NodePotential1} \\
&\pi^i_k \le \pi^i_t + \mu_k, \qquad \qquad && a=(k,t) \in A^i, \label{constr:Level3:Dual2:NodePotential2} \\
&\mu_k  \geq 0, && k \in V, \label{constr:Level3:Dual2:NoneNegMu} \\
&\pi^i_k  \in \mathbb{R}, && k \in V^i. \label{constr:Level3:Dual2:VarDefinition2}
\end{alignat} 
where $\mu_k$ is the dual variable for constraints \eqref{constr:Level3:LocVariables}. As $f_a^i$ is also contained in those, we have to split constraints~\eqref{constr:Level3:Dual:NodePotential} into those for arcs not including $t$, \eqref{constr:Level3:Dual2:NodePotential1}, and those including $t$, \eqref{constr:Level3:Dual2:NodePotential2}, observe that by definition $\ell_{e_a} + \overline \gamma_{e_a} = 0,$ for all $a=(k,t) \in A^i$.

As we can see, since $\overline X$ will become a decision variable, we will end up with a non-linear objective. To linearize the formulation, we first make the following observation about its optimal solutions.

\begin{lemma}
\label{lem:Dual1}
Formulation \eqref{constr:Level3:Dual2:ObjectiveFunction}$-$\eqref{constr:Level3:Dual2:VarDefinition2} has an optimal solution $(\pi,\mu)$ with $\pi_k^i \ge 0$ for $k \in V$, $\pi_t^i=0$, and $\mu_k=0$ whenever $\overline x_k=1$, $k \in V$.
\end{lemma}

\ProofNoNL
The first results follows immediately from the fact that we can replace the equality in \eqref{constr:Level3:FlowBalance} by ``$\ge$'', which results in the $\pi_k^i$-variables being defined as non-negative in the dual.

Let now $(\pi, \mu)$ be an optimal solution of formulation \eqref{constr:Level3:Dual2:ObjectiveFunction}$-$\eqref{constr:Level3:Dual2:VarDefinition2} with $\pi_t^i > 0$ and such that there exists a vertex $j \in V$ with $\overline x_j=1$ and $\mu_j > 0$. 
First, we observe that from the complementary slackness conditions for \eqref{constr:Level3:LocVariables}, $\mu_j(f_{(j,t)}^i - \overline x_j)=0$, we must have $\overline x_j=f_{(j,t)}^i=1$, i.e., the shortest path, $P_{it}$, between $i$ and $t$ must be routed through $j$. Moreover, by construction of $A^i$, $j$ must be the only node of the set $\{k \in V \mid \overline x_k=1 \}$ on $P_{it}$ and with a positive $\mu_k$. Hence, $\sum_{k \in V} \overline x_k \mu_k = \mu_j$.

Second, we observe that from the complementary slackness conditions for \eqref{constr:Level3:Dual2:NodePotential1}, $f_a^i(\pi^i_k - \pi^i_l - \ell_{e_a} - \overline \gamma_{e_a})=0$, we know that for all arcs $a=(k,l)$, $l \neq t$, on $P_{it}$ we must have $\pi^i_k = \pi^i_l + \ell_{e_a} + \overline \gamma_{e_a}$. In addition, from the complementary slackness conditions for \eqref{constr:Level3:Dual2:NodePotential2}, $f_{(j,t)}^i(\pi^i_j - \pi^i_t - \mu_j)=0$, we know that $\pi^i_j = \pi^i_t + \mu_j$.
Hence, $\pi_i^i$ is the length of $P_{it}$ plus $\pi^i_t + \mu_j$, i.e.
\[ \pi_i^i \:=\: \sum_{a=(k,l) \in P_{it}: l \neq t} (\ell_{e_a} + \overline \gamma_{e_a}) + \pi^i_t + \mu_j\,.
\]
With the same argument, $\pi_k^i \ge \pi^i_t + \mu_j$ for all nodes $k \in V$ along $P_{it}$.

Combining these observations, we can construct an alternative optimal solution $(\hat \pi, \hat \mu)$ by setting $\hat \pi_t^i=\hat \mu_j=0$, $\hat \pi_k^i = \pi_k^i - \pi_t^i - \mu_j$, $k \in V$, and $\hat \mu_k=\mu_k$, $k \in V \setminus \{j\}$. 
This solution obviously satisfies \eqref{constr:Level3:Dual2:NodePotential1} and \eqref{constr:Level3:Dual2:NoneNegMu}. It also fulfils \eqref{constr:Level3:Dual2:NodePotential2} for any $k \in V \setminus \{j\}$, as $\hat \pi_k^i=\pi_k^i - \pi_t^i - \mu_j \le \pi_k^i - \pi_t^i \le \mu_k$ and it trivially holds for $j$. Hence, $(\hat \pi, \hat \mu)$ is feasible and, therefore, optimal for \eqref{constr:Level3:Dual2:ObjectiveFunction}$-$\eqref{constr:Level3:Dual2:VarDefinition2} (as $\sum_{k \in V} \overline x_k \mu_k = \mu_j$).
The only thing we need to check is that $\hat \pi_k^i \ge 0$. We already know this holds for all $k \in P_{it}$ from our second observation. Concerning vertices $k \in V$ not on the shortest path, if $\hat \pi_k^i < 0$, we can simply set $\hat \pi_k^i = 0$ without violating any of  \eqref{constr:Level3:Dual2:NodePotential1} (if neither $k$ nor $l$ are on $P_{it}$, then this is trivial, and if, e.g., $l \in P_{it}$, then the right-hand side in \eqref{constr:Level3:Dual2:NodePotential1} is strictly positive) or changing the objective function value, which completes the proof. 
\EndProofNoNL

With this result, we can now state an equivalent formulation that is linear even when $x$ is a decision variable.

\begin{prop}
Formulation \eqref{constr:Level3:Dual2:ObjectiveFunction}$-$\eqref{constr:Level3:Dual2:VarDefinition2} is equivalent to formulation $SP^i$
\begin{alignat}{4}
(SP^i) \qquad d_i(\overline X, \overline \gamma) \:=\: \max\; \pi^i_{i} &  \nonumber\\
\text{s.t. } &\mbox{\labelcref{{constr:Level3:Dual2:NodePotential1}}}, \nonumber\\
\pi^i_k &\le  M (1-\overline x_k), \qquad \qquad && k \in V, \label{constr:Level3:Dual3:NodePotential} \\
\pi^i_k & \geq 0, && k \in V.
\end{alignat} 
\end{prop}
\ProofNoNL
Let $\pi$ be an optimal solution of SP$^i$. Define $(\hat \pi, \hat \mu)$ with $\hat \pi_k^i = \pi_k^i$, $k \in V$, $\hat \pi_t^i=0$, $\hat \mu_k=0$ ($\hat \mu_k = M$) for all $k \in V$ with $\overline x_k=1$ ($\overline x_k=0$). Then, $(\hat \pi, \hat \mu)$ is feasible for formulation \eqref{constr:Level3:Dual2:ObjectiveFunction}$-$\eqref{constr:Level3:Dual2:VarDefinition2} and has the same objective value $\hat \pi_i^i = \pi_i^i$.

In reverse, let $(\pi, \mu)$ be an optimal solution for \eqref{constr:Level3:Dual2:ObjectiveFunction}$-$\eqref{constr:Level3:Dual2:VarDefinition2} that satisfies the properties shown in Lemma~\ref{lem:Dual1}. Then, $\pi$ is feasible for SP$^i$ and has the same objective value.

\EndProofNoNL

Observe that $SP^i$ is formulated on the original network and the use of big-M constraints allows us to avoid the auxiliary network. 

\subsubsection{Middle level problem}
\label{sec:mid_level_problem}

In this subsection, we formulate the attacker decision, named \nameat.  To model the middle level problem, we assume that the upper level decision (location decision) is given and denoted as $\overline X$.  Recall $C(\overline X) = \bigcup_{\overline x \in \overline X} C(\overline x)$ is the set of all nodes that are covered by facilities in $\overline X$ before downgrading. The attacker now aims at finding a downgrade $\gamma$ that maximizes the amount of demand of vertices in $C(\overline X)$ that are no longer covered after the increasing of the length of the edges. In the following, we denote these vertices as being \emph{un-covered} (using the hyphen to avoid ambiguity with the word ``uncovered''). We define the following two families of decision variables:

\begin{table}[h]
\begin{tabular}{lp{0.85\textwidth}}
$\eta_i$ & 1, if node $i \in V$ is un-covered, 0 otherwise, for $i \in V$. \\[0.25ex]
$\rho_i$ & 1, if node $i \in V$ is covered before the downgrading, 0 otherwise, for $i \in V$.\\
\end{tabular}	
\end{table}

\noindent
As we assume the location decisions are given, we set $\overline \rho_i = 1$, for all $i \in C(\overline X)$  and $\overline \rho_i$ = 0 otherwise. The attacker's objective function is to maximize the demand of the un-covered nodes, giving rise to the following formulation:
\begin{alignat}{4}
&& \max \quad& \sum_{i \in V} w_i \eta_i & \label{constr:Level2:Upper:Objective} \\
&& \text{s.t. } \quad& \sum_{e\in E} c_e\gamma_e  \leq B, \label{constr:Level2:Upper:Budget} \\
&&& \eta_i  \leq \overline \rho_i, && i\in V, \label{constr:Level2:Upper:Uncovering} \\
&&& R \eta_i  \leq d_i(\overline X, \gamma), \qquad && i\in V, \label{constr:Level2:Upper:Distance} \\
&&& \eta_i  \in \{0,1\}, && i\in V, \label{constr:Level2:Eta:Definition} \\
&&& 0 \leq \gamma_e  \leq u_e,  && e\in E. \label{constr:Level2:UpperBound:Definition} 
\end{alignat}
The objective function \eqref{constr:Level2:Upper:Objective} maximizes the amount of un-covered demand. Constraints~\eqref{constr:Level2:Upper:Budget} model the budget restriction. Constraints~\eqref{constr:Level2:Upper:Uncovering} and \eqref{constr:Level2:Upper:Distance} ensure that the un-covered nodes were previously covered and that the shortest path distance after the downgrade, i.e., $d^i(\overline X, \gamma)$, is larger than or equal to the coverage radius. 
We point out that adhering to the classical definition of coverage would result in a strict inequality in \eqref{constr:Level2:Upper:Distance}, rendering the feasible region to be non-closed.

Using the formulation $SP^i$ presented in the previous subsection, the attacker's problem can be reformulated as the following single-level problem:
\begin{alignat}{4}
\max \quad &\sum_{i \in V} w_i \eta_i && \nonumber \\
\text{s.t. } &\mbox{\labelcref{constr:Level2:Upper:Budget,constr:Level2:Eta:Definition,constr:Level2:UpperBound:Definition,constr:Level2:Upper:Uncovering}}, \nonumber\\
&R \eta_i  \leq \pi^i_{i},  & & i\in V, \label{constr:Level2:Relation:EtaPi} \\
&\pi^i_k \le \pi^i_l + \ell_{e_a} +  \gamma_{e_a}, \qquad && a=(k,l) \in A, i \in V,\label{constr:Level2:Relation:PiLength} \\
&\pi^i_k \le  M (1-\overline x_k), \qquad && i, k \in V, \label{constr:Level2:Relation:BigM} \\
&\pi^i_k  \geq 0, && i, k\in V. \label{constr:Level2:Pi:Definition}
\end{alignat} 

If each node $i \in V$ has a unique closest facility in $\overline X$, then the set of shortest paths of each node to its closest facility naturally defines a forest, with each tree of the forest rooted at a facility in $\overline X$. If a node $i \in V$ is equidistant to two or more facilities in $\overline X$, then we again obtain a forest by assigning $i$ without loss of generality to one of those facilities as well as all other nodes for which this facility is among their closest.
As a result, the potentials of nodes in each tree can be computed independently from one another. Moreover, within each tree the potential of a node $k$ can be the same for all nodes in the tree, i.e., we can assume \wloge that $\pi_k^i=\pi_k^j=\pi_k$ with $i$ and $j$ being nodes of the tree, as the potential of a node can represent the length of the path from $k$ to the root of the tree. Therefore, we do not need to define the $\pi$-variables for each $i\in V$, i.e., we can remove the super-index $i$ resulting in the reduced formulation:
\begin{alignat}{4}
Q(\overline x, \overline \rho) = \max \quad &\sum_{i \in V} w_i \eta_i && \nonumber \\
\text{s.t. } &\mbox{\labelcref{constr:Level2:Upper:Budget,constr:Level2:Eta:Definition,constr:Level2:UpperBound:Definition,constr:Level2:Upper:Uncovering}}, \nonumber\\
&R \eta_i  \leq \pi_{i},  & & i\in V, \label{constr:Level2:Relation:EtaPi_1ind} \\
&\pi_k \le \pi_l + \ell_{e_a} +  \gamma_{e_a}, \qquad && a=(k,l) \in A, \label{constr:Level2:Relation:PiLength_1ind} \\
&\pi_k \le  M (1-\overline x_k), \qquad && k \in V, \label{constr:Level2:Relation:BigM_1ind} \\
&\pi_k  \geq 0, && k\in V. \label{constr:Level2:Pi:Definition_1ind}
\end{alignat}

The previous formulation is enhanced in the next subsection, and a valid value for the big-M is provided. In addition, further improvements are discussed in Section~\ref{sec:Heuristic}. However, they take advantage that the location of the facilities is known, so they can not be applied to the general problem. 

\subsubsection{Upper level problem}

In this subsection, we introduce the formulation of \name called \namef. It is a linear bilevel formulation, where the objective function is to maximise the amount of covered demand after the attacker has optimally downgraded the network. 
\begin{alignat}{4}
\max \, &\sum_{i\in V} w_i \rho_i  - Q(x,\rho), \label{constr:Level1:Upper:Objective_1} \\
\text{s.t. } &\sum_{j\in V} x_j  = p, && \label{constr:Level1:Upper:NbFacilities_1} \\
&\rho_i  \leq \sum_{j\in V:d(i,j)<R} x_j, \qquad && i\in V, \label{constr:Level1:Upper:Coverage_1} \\
&x_i \in \{0,1\}, && i\in V,\\
&\rho_i  \in \{0,1\}, && i\in V,\\
&Q(x,\rho)  = \max \sum_{i \in V} w_i \eta_i \label{constr:Level1:Upper:Attacker_1} \\
& \qquad \text{s.t. } \mbox{\labelcref{{constr:Level2:Upper:Budget,constr:Level2:Eta:Definition,constr:Level2:UpperBound:Definition,constr:Level2:Pi:Definition_1ind,constr:Level2:Relation:PiLength_1ind,constr:Level2:Relation:EtaPi_1ind}}}, \nonumber\\
&\qquad\qquad\eta_i  \leq  \rho_i, && i\in V,\label{constr:Level1:Lower:Defender_1}\\
&\qquad\qquad \pi_k \le  M_k (1- x_k), \qquad && k \in V. \label{constr:Level1:Lower:Defender_2_1ind} 
\end{alignat} 

The objective function \eqref{constr:Level1:Upper:Objective_1} maximizes the amount of demand that is still covered after the attacker acted. Constraints~\eqref{constr:Level1:Upper:NbFacilities_1} and \eqref{constr:Level1:Upper:Coverage_1} ensure that exactly $p$ facilities are located and, respectively, that a node $i$ can only be initially covered if there is a facility within the coverage radius. Constraint~\eqref{constr:Level1:Upper:Attacker_1} calculates the amount of demand un-covered by the attacker. This is also the objective function of the follower, whose associated constraints were described in the previous section.  Furthermore, Constraints~\eqref{constr:Level1:Lower:Defender_1} establishes the relationship between the leader variable $\rho$ and the follower variables $\eta$. Similarly, Constraints~\eqref{constr:Level1:Lower:Defender_2_1ind} sets the relationship between the leader variable $x$ and the lower level variables $\pi$, whose associated big-M depends on index $k$, to make it tighter.

\begin{lemma}
A valid value for $M_k$ in constraint~\eqref{constr:Level1:Lower:Defender_2_1ind} is the 
$(n-p+1)$-th distance $d_u(k,j)$ sorted in non-decreasing order for $j\in V$, where $d_u(k,j)$ is the distance from node $k$ to node $j$ in a network $N(u)=(V, E, \ell(u))$.
\end{lemma}
\ProofNoNL
A valid value for the big-M in constraint~\eqref{constr:Level1:Lower:Defender_2_1ind} is an upper bound of variable $\pi_k.$ Recall that this variable represents the node potential of node $k$. As we fix the node facilities' potential to zero, we can assume \wloge that the potential of a node $k$ is the distance from $k$ to the closest service facility. As, a priori, the node that will be the facility is not known, the length from $k$ to the furthest node is a valid upper bound of the potential of a node $k$. Indeed, if $p$ facilities are located, we can assume that the distance from a node to the closest facility is at most the distance to the $(n-p+1)$-th furthest node.  Similarly, the downgrading to be applied is unknown, so a worst-case scenario is assumed, i.e., for computing the distance we suppose that all edges have been downgraded to their maximum ($\ell_e(u)=\ell_e+u_e,$ for $e\in E$). Therefore, we obtain that a valid value for $M_k$ is the $(n-p+1)$-th distance $d_u(k,j),$ sorted in non-decreasing order for $j\in V$.
\EndProofNoNL

\section{Preprocessing} \label{sec:Prep}

In this section, we present a strategy for preprocessing data to reduce the number of variables and constraints included in the model. The procedure is based on obtaining transformed input data, where the optimal solution is also optimal in the original data. To achieve this, we follow the steps outlined below:

\begin{description}
    \item [Identification of redundant edges:] When solving the $\name$, any edges $e\in E$ in the original network with lengths exceeding $R$ can be safely removed without impacting the optimal solution. Consequently, their associated costs and upper bounds can also be removed, i.e., $c_e$ and $u_e$. 
    \item [Tightening downgrading upper bounds ($u$):] The attacker's objective is to compromise the leader's coverage. Thus, when the distance from the facility to the client equals $R$, the attacker achieves its goal. Consequently, if there exists an edge $e\in E$ whose length plus the upper bound exceeds $R$, i.e., $\ell_e+u_e\geq R$ the upper bound can be tightened as follows:
    \begin{equation*}
        u_e=R-\ell_e.
    \end{equation*}
 
\end{description}

Therefore, removing all these edges eliminates numerous constraints from the formulation. Additionally, reducing the upper bounds of the $\gamma$ variables improves the bounds of the formulation. In the following, unless stated otherwise, we refer to \namef as the bilevel formulation where these enhancements have been applied. 

As reported in Section~\ref{sec:ComputationalResults}, this strategy is useful for solving exactly our problem on the general purpose solver for bilevel programming \citep{FisLjuMonSin17}, since it notably reduces the number of constraints. Furthermore, it also benefits the matheuristic that is introduced in the following section. 

\section{Matheuristics}
\label{sec:Heuristic}

The proposed heuristic approach is based on generating feasible solutions. In this sense, it is worth mentioning that to provide
a feasible solution of \name, the attacker problem A-MCLP must be solved optimally for a given set \ox of facility locations. Otherwise, the solution may not be feasible. 
However, as the next lemma shows, finding the optimal edge length increases for a given \ox is NP-hard.
\begin{lemma}
\label{lemma:NPhard:AttackerProblem}
Given a set \ox of facility locations, solving \nameat is NP-hard even on star networks for non-uniform weights assuming uniformity in downgrading costs and upper bounds as well as integrality of the input parameter values.
\end{lemma}

\noindent
\ProofNoNL See Appendix~\ref{sec:App}.

\medskip

Therefore, we first focus on improving the formulation of the attacker problem. Let a set \ox of facility locations be given. Then, we define the set $\overline{V}:=\{i\in V: \overline \rho_i=1\},$ i.e., the covered nodes in the original network. Note that $\overline \rho_i$ can be computed easily for a given \ox, e.g. using Dijkstra's algorithm starting from each $x \in \overline{X}$. Below, for any $m\in V$, we define some sets that can be used to reduce the number of constraints of the formulation.
\begin{align*}
      &V^m:=\{i\in V: d(i,m)< R\}, \\
  &A[V^m]:=\{(k,l)\in A: k,l \in V^m\}.
\end{align*}
The first set includes the nodes that can be covered by $m$ (before downgrading) while the second comprises the arcs that can be used for this covering. For the ease of exposition, we rewrite the formulation for the middle level problem:
\begin{alignat}{4}
\max \quad \sum_{i \in \overline{V}} w_i & \eta_i && \nonumber \\
\text{s.t. } &\mbox{\labelcref{{constr:Level2:Upper:Budget,constr:Level2:UpperBound:Definition}}}, \nonumber\\ 
& R \eta_i  \leq \pi_{i},  & & i\in \overline V, \label{cons:matheuristic:eta_pi} \\
&\pi_k \le \pi_l + \ell_{e_a} +  \gamma_{e_a}, \qquad && k,l \in \overline V: a=(k,l) \in \bigcup_{m\in \overline X} A[V^m], \label{cons:relatespi_2}\\
&\pi_k =0,\qquad && k \in \overline X, \label{cons:matheuristic:BigM} \\
&\pi_k  \geq 0, && k \in \overline V \setminus \overline X, \label{cons:matheuristic:BigM2} \\
&\eta_i  \in \{0,1\},  && i\in \overline V. \label{cons:matheuristic:Definition_eta}
\end{alignat}

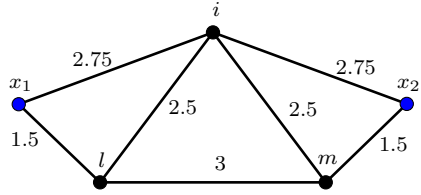
\begin{figure}[htp]  \label{DG:example}
\centering
\begin{tikzpicture}[line cap=round,line join=round,x=1.0cm,y=1.0cm]
\clip(2.294526104206216,-0.2545734346927407) rectangle (9.167479946523589,2.5107194290970316);
\draw [line width=1.pt] (4.,0.)-- (7.,0.);
\draw [line width=1.pt] (4.,0.)-- (5.50095073444791,1.9949635522469436);
\draw [line width=1.pt] (5.50095073444791,1.9949635522469436)-- (7.,0.);
\draw [line width=1.pt] (5.50095073444791,1.9949635522469436)-- (8.07713946451166,1.0439207699780355);
\draw [line width=1.pt] (7.,0.)-- (8.07713946451166,1.0439207699780355);
\draw [line width=1.pt] (5.50095073444791,1.9949635522469436)-- (2.9184838974351415,1.0393858378354586);
\draw [line width=1.pt] (4.,0.)-- (2.9184838974351415,1.0393858378354586);
\begin{scriptsize}
\draw [fill=black] (4.,0.) circle (2.5pt);
\draw[color=black] (4.0273083002439085,0.3) node {$l$};
\draw [fill=black] (7.,0.) circle (2.5pt);
\draw[color=black] (7.027473157652096,0.3) node {$m$};
\draw[color=black] (5.602187082733803,0.25) node {$3$};
\draw [fill=black] (5.50095073444791,1.9949635522469436) circle (2.5pt);
\draw[color=black] (5.1,0.996093264934593) node {$2.5$};
\draw[color=black] (6.7,0.9753164999940932) node {$2.5$};
\draw[color=black] (5.544012140900403,2.3) node {$i$};
\draw [fill=blue] (8.07713946451166,1.0439207699780355) circle (2.5pt);
\draw[color=black] (8.107864934558092,1.3) node {$x_2$};
\draw[color=black] (7.4,1.5321338003994889) node {$2.75$};
\draw[color=black] (7.9,0.5099169653268968) node {$1.5$};
\draw [fill=blue] (2.9184838974351415,1.0393858378354586) circle (2.5pt);
\draw[color=black] (2.9469165233379133,1.3) node {$x_1$};
\draw[color=black] (3.9,1.6360176251019882) node {$2.75$};
\draw[color=black] (3,0.5805579661245962) node {$1.5$};
\end{scriptsize}
\end{tikzpicture}
\caption{Graphical representation of the improvement.}
\end{figure}

Constraints \eqref{cons:matheuristic:eta_pi}, \eqref{cons:matheuristic:BigM}, \eqref{cons:matheuristic:BigM2}, and \eqref{cons:matheuristic:Definition_eta} replace \eqref{constr:Level2:Relation:EtaPi_1ind}, \eqref{constr:Level2:Relation:BigM_1ind}, \eqref{constr:Level2:Pi:Definition_1ind}, and \eqref{constr:Level2:Eta:Definition}, respectively, and \eqref{constr:Level2:Upper:Uncovering} becomes redundant due to replacing $V$ by $\overline V$. Moreover, \eqref{constr:Level2:Relation:PiLength_1ind} may be restricted to those arcs whose end nodes can be covered by the same open facility $m \in \overline X$ before the downgrade, resulting in \eqref{cons:relatespi_2}. We can do this because an arc for which this does not hold can not be included in a path from any node to its covering facility with a length shorter than $R$. Note that we assumed that the redundant edges (the edges such that $\ell_e\geq R)$ have already been removed from the network. Therefore, this represents a notable reduction in the number of constraints included in the formulation. In Figure~\ref{DG:example}, a graphical example is included. For $R$ equal to $4$, if the facility is located at $x_1$ (respectively at $x_2$), then $i$ and $l$ will be covered by $x_1$ ($i$ and $m$ by $x_2$). Note that the edge $(l,m)$ represented in Figure \ref{DG:example} will not be included in the constraint, as the same facility cannot cover both end nodes and, thus, it can never be part of a shortest path from $i$ to $x_1$ (respectively $x_2$) whose length is $<R$. 

Before we proceed to describe the phases of our matheuristic, we first note that we can easily compute a lower bound and an upper bound of \name. Observe that these bounds may be unattainable, as they may not be associated with a feasible solution. To do that, we will make extensive use of auxiliary networks where edge lengths have been modified according to given values. Let $\gamma \in \mathbb{R}^{m}$ such that $0 \le \gamma_e \le u_e$, $e \in E$. Recall $N(\gamma)=(V,E,\ell(\gamma))$ denotes a network with edge lengths $\ell_e(\gamma)=\ell_e+\gamma_e$.

\begin{description}
\item[Lower bound:] We obtain a lower bound by solving the classical MCLP in the auxiliary network 
$N(u)$. The optimal objective value of this auxiliary problem is a lower bound of \name,  as it is clear that the leader can cover at least this amount of demand. We call the corresponding set of facility locations $\overline X_{LB}.$ Observe that this downgrading strategy may not be feasible for the attacker due to the limited downgrading budget $B$. Basically, this lower bound can be obtained by solving the attacker problem assuming that $B=\displaystyle \sum_{e\in E} c_eu_e.$ 

\item[Upper bound:] We obtain an upper bound by solving the classical MCLP in the original network, i.e., we assume no edges have been downgraded ($\gamma_e=0, e\in E$). The optimal objective value is an upper bound \name, as it is evident that the leader can cover at most this demand. We call the corresponding set of facility locations $\overline X_{UB}.$ Note that this downgrading strategy will likely not be optimal for the attacker for values of $B$ larger than 0. Basically, this upper bound can be obtained by solving the attacker problem with $B=0.$ 
\end{description}

\noindent
Our heuristic is split into two phases: a construction phase (alternating location-downgrading search) and an improvement phase (1-1 local search).

\subsection{Alternating Location-Downgrading Search}
In this subsection, we describe the first phase of the algorithm: the alternating location-downgrading search. This is a repetitive approach where we continuously build feasible solutions in a series of iterations until no further improvements can be achieved or a maximum number of iterations have been reached. The steps are described below:
    
\begin{description}
    \item[Step 0] Propose a promising set $\overline X$ of facility locations.

    \item[Step 1] 
    Solve the attacker problem A-MCLP for given locations $\overline X$ to get an optimal downgrade $\overline \gamma$. $\overline X$ and $\overline \gamma$ constitute a feasible solution of $\name$.  

    \item[Step 2] Solve the classical MCLP in $N(\og)$ to obtain an optimal set of locations \ox.

    \item[Step 3] If \ox has not changed from the previous iteration or the maximum number of iterations has been reached, then stop; otherwise, return to Step 1.
\end{description}

We start this alternating location-downgrading search from the following promising sets of facility locations (we depict the ``code" of the strategy in brackets). We choose $\overline X$ as an optimal solution of the classical MCLP in the network $N(\gamma^0)$, where $\gamma^0$ is: (0) $\gamma_e^0=0$, $e \in E$, i.e., $\overline X=$ $\overline X_{LB}$, (1) $\gamma_e^0=u_e,$ i.e., $\overline X=$ $\overline X _{UB},$ (2)  $\gamma_e^0=\min\left\{\frac{B}{|E|}, u_e\right\},$ i.e., proportional to the budget and inversely proportional to the number of edges, (3) $\gamma_e^0=\min\left\{\frac{u_e B}{\sum_{e\in E }u_e}, u_e\right\},$ i.e., proportional to the budget and the upper bounds, (4) proportional to the benefit-cost ratio, i.e., downgrade the cheapest edges first until reaching the budget, (5,6,7) solve the classical MCLP considering $80\%R, 70\%R, 60\%R$, respectively, and (8) $\gamma_e^0=u_e/2, e\in E,$ i.e., all the edges have been downgraded to half of their maximum.

We run the alternating location-downgrading search with each of the eight starting solutions and we keep the best set of locations \ox and their corresponding optimal downgrades \og.

\subsection{1-1 Local Search}

Once the previous iterative process has finished with a feasible solution \ox and \og, we developed a 1-1 local search to improve this solution by substituting one facility location in $\overline X$. For doing so, we look for a node $i \in \overline X$ to be removed as facility and we compute which is the most promising node $j \not\in \overline X$ to be included. This local search is repeated as long as better feasible solutions are found, or up to a fixed number of iterations. Unfortunately, optimally solving the attacker problem numerous times can be very costly. Therefore, we propose three different strategies for evaluating a new location set, so the user can select the best strategy according to the trade-off between computation time and accuracy required. Moreover, each strategy includes several variants, resulting in a diverse set of matheuristics. The faster, but less refined strategies are detailed first.

\begin{description}
\item[Fixed Out-In] $ $

First, to select the node $i \in \overline X$ to leave the facility location set, we consider the network where the edge lengths have been modified according to the optimal downgrading strategy \og for the current set \ox.  

In the modified network $N(\og)$ we compute the amount of covered demand considering the removal of each facility location. Then, we choose as the most promising node $j \not\in \ox$ to be included in the set $\ox \setminus \{i\}$, for each $i \in \ox$, the node whose amount of covered demand is the largest in the following two networks (also considering the demand covered in $N(\og)$ by the nodes remaining in the facility location set): 
\begin{enumerate}
  \item[a)] $N(\og)$, where $\og$ is the optimal downgrading strategy for the current  location set \ox. Observe that the covered demand is computed taking into account simultaneously the node that leaves and the one that enters the location set.
  \item[b)] $N(u)$. Observe that the demand covered by the nodes remaining in the facility location set is computed in $N(\og)$. Then, the covered demand is computed by simultaneously considering the node that leaves and the node that enters the location set, each in a different network.

\end{enumerate} 

Finally, the most beneficial output-input pair is selected. For this new set of facilities, the \nameat problem is solved to obtain a feasible solution of \name. This procedure is summarized in Algorithm~\ref{algo:LS_FOI}.

\begin{algorithm}
\SetKwComment{Comment}{/* }{ */}
\KwData{Feasible solution $(\overline X, \overline \gamma)$ of \name.}
\KwResult{New set \ox of candidate locations with optimal downgrades \og.}
\BlankLine
\For{$i\in \overline{X},j \in V \setminus \overline{X}$}{
 
    a) Compute $\sum_{k\in \mathcal{C}\left((\overline X\setminus \{i\})\cup \{j\}, \overline \gamma\right)} w_k.$ \Comment*[r]{For version a).}
    
    b) Compute $\sum_{k\in \mathcal{C}(\overline X\setminus \{i\}, \overline \gamma) \cup \, \mathcal{C}\left(\{j\}, u\right)}w_k.$ \Comment*[r]{For version b).}
  }

Select the pair $i,j$ that provides the largest amount of coverage. 

Set $\overline X=(\overline X\setminus \{i\}) \cup \{j\}$ and compute the optimal downgrades \og for \ox.
\caption{Local Search: Fixed Out-In, version a) and b).}\label{algo:LS_FOI}
\end{algorithm}

\item [Fixed Out-Optimal In] $ $
    
Again, we have to choose one node to leave and one node to enter the location set. 
In contrast to the previous procedure, this approach first fixes the leaving node $i \in \ox$ and once it is fixed, the most promising node $j \not\in \ox$ to enter is obtained. 
To find the candidate to be removed, we select the node $i^* \in \ox$ whose removal reduces the amount of covered demand by the least in the following two networks:
\begin{enumerate}
  \item[a)] $N(\og)$. 

  \item[b)] $N(u)$.  
 
\end{enumerate} 

Now, we look for the candidate to enter the location set $\ox \setminus \{i^*\}$. In this case, we test all nodes $j \not\in \ox$. For each node, we optimally solve the attacker problem for the set of facilities where the selected node in the strategy a) or b) is removed and the candidate node is included, i.e., for $(\ox \setminus \{i^*\}) \cup \{j\}$. Note that if the attacker problem for this solution set has already been solved, it is not tested again. Finally, we select the candidate that provides the best objective function value of \name.

A variant of this strategy is solving the attacker problem for a given location of the facilities with a time limit, e.g., 10 seconds. A feasible solution to the attacker problem will overestimate the actual objective function value of \name for the set $(\ox \setminus \{i^*\}) \cup \{j\}$ but it can notably save computational time. Once a candidate $j^*$ to join the set of locations has been chosen, the attacker problem is solved exactly for these locations. 

The procedure is summarised in Algorithm~\ref{algo:LS_FOoptI}.
\begin{algorithm}
\SetKwComment{Comment}{/* }{ */}
\KwData{Feasible solution $(\overline X, \overline \gamma)$ of \name.}
\KwResult{New set \ox of candidate locations with optimal downgrades \og.}
\BlankLine
\For{$i\in \overline X$}{
    a) Compute $\sum_{k\in \mathcal{C}\left(\{i\}, \overline \gamma\right) \setminus \mathcal{C}\left(\overline X\setminus \{i\}, \overline \gamma\right)} w_k.$ \Comment*[r]{For version a).}
    
    b) Compute $\sum_{k\in \mathcal{C}\left(\{i\}, u\right) \setminus \mathcal{C}\left(\overline X\setminus \{i\}, u\right)}w_k.$ \Comment*[r]{For version b).}
}

Select the node $i^*$ that provides the previous smallest amount to leave the location set. 

\For{$j \in V \setminus \overline X$}{

    Compute the optimal downgrades \og for the set of locations $(\overline X\setminus \{i^*\}) \cup \{j\}$ and the objective value of \name.

}
Select the node $j^*$ that provides the largest objective value to enter the location set. 
  
Set $\overline X=(\overline X\setminus \{i^*\}) \cup \{j^*\}.$ 
\caption{Local Search: Fixed Out-Optimal In, version a) and b).}\label{algo:LS_FOoptI}
\end{algorithm}

\item [Optimal Out-In] $ $

Finally, we present the last strategy, which requires the most computational resources. In this case, we consider all pairs $(i,j)$ of one candidate $i$ leaving the set \ox of locations and one candidate $j \not\in \ox$ entering the set. For each pair, we optimally solve the attacker problem for the new location set $(\ox \setminus \{i\})\cup\{j\}$. Finally, we select the candidate to leave and the candidate to enter that provides the best objective function value of \name. 

As before, a variant of this strategy consists of solving the attacker problem with a time limit and only solving the attacker problem exactly when the candidates to leave and enter the location set has been chosen. 

This procedure is summarised in Algorithm~\ref{algo:LS_OptOI}.
\begin{algorithm}
\SetKwComment{Comment}{/* }{ */}
\KwData{Feasible solution $(\overline X, \overline \gamma)$ of \name.}
\KwResult{New set \ox of candidate locations with optimal downgrades \og.}
\BlankLine
\For{$i\in \overline{X},j \in V\setminus \overline X$}{
 
    Compute the optimal downgrades \og for the set of locations $(\overline X\setminus \{i\}) \cup \{j\}$ and the objective value of \name.

  }
Select the pair of nodes $i$ and $j$ that provides the largest objective value of \name to leave and enter the location set, respectively. 
  
Set $\overline X=(\overline X\setminus \{i\}) \cup \{j\}$. 
\caption{Local Search: Optimal Out-In.}\label{algo:LS_OptOI}
\end{algorithm}
\end{description}

\section{Computational Results} \label{sec:ComputationalResults}

This section is dedicated to computational experiments. These are divided into three main groups, the exact resolution of the formulation using the general solver developed in \cite{FisLjuMonSin17}, the matheuristic experiments, and the managerial insight.

\subsection{Data}\label{sub:data}

We generated  instances adapting the procedure used in
\cite{ReSchoWi08, CorFurLju19, BalKalMarRod22}, among others. The nodes were given by points whose coordinates followed a uniform distribution over [0,30]. We then built the complete graph, where the length of the edges is the Euclidean distance between the nodes rounded to two decimal places.  These instances are called ``graph" followed by $n$ (the number of vertices); for example, ``graph50" is a complete graph with 50 nodes and 1225 edges. 

The selection of parameters is outlined as follows. The number of facilities denoted as $p$ was determined in proportion to the number of vertices, specifically, $p\in \{n/30, n/20, n/10\}$. The weights or demands assigned to nodes, denoted as $w_i$ for $i\in V$, were uniformly randomly generated integers from 1 to 100. We tested three different coverage radii, $R$, such that each node can cover at least one node, at least $5\%$ of the number of nodes, and at least $10\%$ of the number of nodes, respectively. Downgrading costs, $c_e,$ for $e\in E,$ were uniformly randomly generated between 1 and 3 with two decimal places. The upper bounds $u_e$, for $e\in E,$ were uniformly randomly generated from $(0.5\ell_e, 1.5\ell_e),$ for $e\in E$ with two decimal places. 

Finally, the budget $B$ was computed as follows. We calculated the maximum  required  budget for downgrading all the edges, $B_{max}=\sum_{e\in E} u_ec_e,$ selected $B_{per}\in \{0.1, 0.05, 0.025\},$ and computed $B$ as: 

$$B=\frac{B_{max}\cdot B_{per} \cdot p \cdot (p-1)}{n (n-1)} $$ rounded to two decimal places. 

For each combination of parameters, five instances were generated with the same procedure, varying only the random seed for the generator. The instances are available at \cite{Data}. 

\subsection{Results}
In this subsection, we analyse the result provided by the exact formulation, the matheuristic, and the managerial insight of the model.

\subsubsection{Exact formulation}
As stated before, the formulation \namef was tested in the general bilevel solver developed in \cite{FisLjuMonSin17}. The computational experiments were carried out on a computer with an Intel(R) Xeon(R) CPU E5-1650 v4 3.60GHz processor and 16 GB RAM. 

In Table \ref{tab:Exactn_50_p}, formulation \namef  without and with preprocessing is compared for 50 nodes with a time limit of one hour. Moreover, the results of 10 hours of computation using the preprocessing are also reported. The first column of the table reports the number of facilities ($p$). The first group of columns reports the results without using the preprocesing, where the first column indicates the number of instances solved to optimality within the time limit of one hour (\#Sol.). Note that 45 instances were analysed for each $p,$ setting three different values for the budget, three for the covering radius, and five random seeds. The next column indicates the average solution time in seconds “t(s.)”. Observe that for the instances where optimality was not proven within the time limit, a solution time equal to the time limit was assigned to compute the average solution time, i.e., the actual solution times are underestimated. Finally, the average relative percentage gap as reported by the solver is depicted in column (\%G). The following two groups of columns report the same information, but for the formulation with preprocessing and having an hour and 10 hours of time limit, respectively. 
\begin{table}[htbp]
      \caption{Comparison of \namef for $n=50$ by $p$, with and without preprocessing.}
\medskip
        \centering
    \begin{tabular}{l|rrr|rrr|rrr}
          & \multicolumn{3}{c|}{Without prep. (1h)}& \multicolumn{3}{c|}{With prep. (1h)}  & \multicolumn{3}{c}{With prep. (10h)} \\
\cline{2-10}    $p$     & \multicolumn{1}{r}{\#Sol. } & \multicolumn{1}{r}{t(s.)} & \multicolumn{1}{r|}{\%G} & \multicolumn{1}{r}{\#Sol. } & \multicolumn{1}{r}{t(s.)} & \multicolumn{1}{r|}{\%G} &\multicolumn{1}{r}{\#Sol. } & \multicolumn{1}{r}{t(s.)} & \multicolumn{1}{r}{\%G} \\
    \hline
    2   &   45    & 251.70 & 0.00  & 45    & 141.33 & 0.00  & 45    & 141.33 & 0.00   \\
   3 & 15    & 2564.16 & 9.53  & 18    & 2452.20 & 8.16  & 29    & 14997.76 & 4.78 \\
    5    & 5     & 3282.95 & 26.31 & 5     & 3265.33 & 25.25 & 5 & 32058.86& 22.52 \\
    \hline
    \end{tabular}%
  \label{tab:Exactn_50_p}%
\end{table}%
Two similar tables, \ref{tab:Exactn_50_R} and \ref{tab:Exactn_50_B}, are reported below, where instead of grouping instances by the number of services to be located, they are grouped by the covering radius and budget, respectively. Observe that, as before, each line represents the average over 45 instances.

\begin{table}[htbp]
  \centering
        \caption{Comparison of \namef for $n=50$ by $R$, with and without preprocessing.}

        \medskip 
        
  	\resizebox{\textwidth}{!}{
    \begin{tabular}{l|rrr|rrr|rrr}
         & \multicolumn{3}{c|}{Without prep. (1h)}  & \multicolumn{3}{c|}{With prep. (1h)} & \multicolumn{3}{c}{With prep. (10h)} \\
\cline{2-10}    $R$     & \multicolumn{1}{r}{\#Sol. } & \multicolumn{1}{r}{t(s.)} & \multicolumn{1}{r|}{\%G} & \multicolumn{1}{r}{\#Sol. } & \multicolumn{1}{r}{t(s.)} & \multicolumn{1}{r|}{\%G} &\multicolumn{1}{r}{\#Sol. } & \multicolumn{1}{r}{t(s.)} & \multicolumn{1}{r}{\%G} \\
\cline{1-10}    
    At least one     & 19    & 2299.56 & 22.50 & 19    & 2170.07 & 20.60  &24	&17661.29 &16.40 \\
    At least 5\%    & 19    & 2157.99 & 10.42 & 22    & 2088.59 & 10.04 & 24	&17433.81	&8.64 \\
     At least 10\%     & 27    & 1641.26 & 2.91  & 27    & 1600.21 & 2.78  & 31	&12093.17 &	2.27\\
    \hline
    \end{tabular}%
    }

  \label{tab:Exactn_50_R}%
\end{table}%

\begin{table}[htbp]
  \centering
    \caption{Comparison of \namef for $n=50$ by $B_{per}$, with and without preprocessing.}
    \medskip
    
  \begin{tabular}{l|rrr|rrr|rrr}
         & \multicolumn{3}{c|}{Without prep. (1h)} & \multicolumn{3}{c|}{With prep. (1h)}  & \multicolumn{3}{c}{With prep. (10h)} \\
\cline{2-10}    $B_{per}$     & \multicolumn{1}{r}{\#Sol. } & \multicolumn{1}{r}{t(s.)} & \multicolumn{1}{r|}{\%G} & \multicolumn{1}{r}{\#Sol. } & \multicolumn{1}{r}{t(s.)} & \multicolumn{1}{r|}{\%G} &\multicolumn{1}{r}{\#Sol. } & \multicolumn{1}{r}{t(s.)} & \multicolumn{1}{r}{\%G} \\
\cline{1-10}    0.025 & 28    & 1431.76 & 5.00  & 30    & 1378.18 & 4.62  & 33	&10304.86	&3.60 \\
    0.050 & 22    & 2037.98 & 10.67 & 23    & 1963.60 & 10.07 & 27	& 15318.54	&8.17 \\
    0.100 & 15    & 2629.07 & 20.17 & 15    & 2517.09 & 18.72 & 19	& 21564.87 &	15.54 \\

    \hline
    \end{tabular}%

  \label{tab:Exactn_50_B}%
\end{table}%

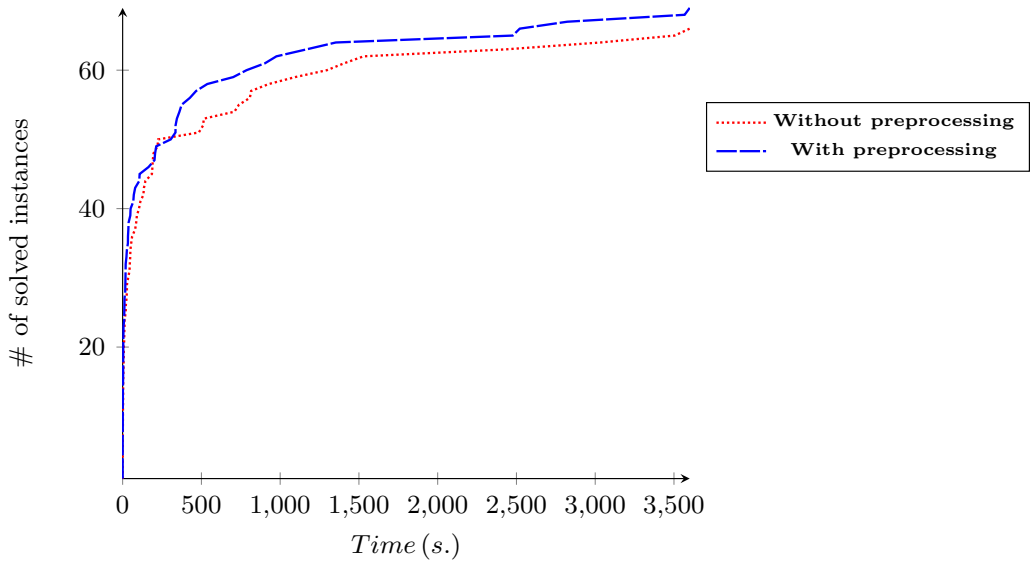
\begin{figure}[hbt] \centering
	\begin{tikzpicture}[scale=1.1,font=\footnotesize]
		\begin{axis}[axis x line=bottom,  axis y line=left,
			xlabel=$Time\,(s.)$,
			ylabel=\# of solved instances,
			legend style={at={(1.6,0.8)}}]
			\addplot[red,thick, densely dotted] plot coordinates {		
(0.084624,1)
(0.108124,2)
(0.187603,3)
(0.412613,4)
(0.413127,5)
(0.447631,6)
(0.718837,7)
(1.237387,8)
(1.648269,9)
(1.748804,10)
(1.929616,11)
(2.873537,12)
(3.104238,13)
(3.268095,14)
(3.505281,15)
(4.1975,16)
(5.852668,17)
(5.900977,18)
(7.34118,19)
(7.529926,20)
(10.045428,21)
(10.330351,22)
(11.917428,23)
(12.238241,24)
(19.665042,25)
(20.084662,26)
(24.405295,27)
(27.08523,28)
(27.211711,29)
(34.980331,30)
(43.563482,31)
(45.510039,32)
(48.030638,33)
(49.061308,34)
(54.265779,35)
(60.325187,36)
(79.460691,37)
(86.035803,38)
(90.466209,39)
(103.7588,40)
(112.311582,41)
(131.673058,42)
(135.235409,43)
(144.6744,44)
(185.733372,45)
(189.111009,46)
(191.264015,47)
(196.433185,48)
(221.540607,49)
(225.205882,50)
(483.475093,51)
(510.410088,52)
(515.622995,53)
(709.529768,54)
(737.943017,55)
(807.701128,56)
(815.44475,57)
(928.008408,58)
(1095.235634,59)
(1296.832546,60)
(1406.698045,61)
(1529.539254,62)
(2427.412791,63)
(3032.669616,64)
(3504.327975,65)
(3600.004738,66)
			};
			\addlegendentry{\tiny\textbf{Without preprocessing}}	
   \addplot[blue,thick,dash pattern=on 7pt off 1pt] plot coordinates {		
(0.037081,1)
(0.051351,2)
(0.114925,3)
(0.257168,4)
(0.27462,5)
(0.287112,6)
(0.341032,7)
(0.61216,8)
(0.901355,9)
(0.925271,10)
(0.945848,11)
(1.060839,12)
(1.532338,13)
(1.687327,14)
(2.516037,15)
(3.086909,16)
(3.208093,17)
(3.695035,18)
(3.810021,19)
(4.015866,20)
(4.204751,21)
(5.095502,22)
(5.433484,23)
(10.345544,24)
(10.419314,25)
(11.423533,26)
(11.88382,27)
(16.451156,28)
(16.961574,29)
(18.010625,30)
(18.683763,31)
(19.28824,32)
(24.021282,33)
(26.75446,34)
(33.484368,35)
(33.8959,36)
(36.314653,37)
(37.916557,38)
(48.60514,39)
(50.352722,40)
(67.445763,41)
(70.834144,42)
(79.436473,43)
(105.472352,44)
(107.891183,45)
(164.073491,46)
(202.14658,47)
(206.041828,48)
(213.291115,49)
(305.395477,50)
(335.449179,51)
(335.824208,52)
(344.343841,53)
(360.635721,54)
(374.309257,55)
(426.79122,56)
(466.548854,57)
(537.315114,58)
(702.09576,59)
(786.785965,60)
(903.735164,61)
(976.269758,62)
(1165.99629,63)
(1353.212848,64)
(2476.575021,65)
(2520.949004,66)
(2822.654964,67)
(3567.617262,68)
(3600.003857,69)

			};
			\addlegendentry{\tiny\textbf{With preprocessing}}
  
		\end{axis}
	\end{tikzpicture}
	\caption{Performance profile graph of \#solved instances (out of 135) for $n=50$.}  \label{fig:performanceProfilesALL}
\end{figure}

The results show the computational efficiency of the preprocessing. Theoretically, it was already shown that it reduced the number of variables and constraints of the formulation. Now, it can be observed that using the preprocessing enables us to solve more instances to optimality, reducing the average computational time, as well as reducing the GAP of the solution. For example, for $p=3,$ the formulation without preprocessing solved to optimality within the time limit 15 instances while using the preprocessing 18 instances were solved. The benefits of preprocessing are evident in all the values of the number of services to be located, and also when grouped according to different radii and budgets. Overall, it is found that about 48.15\% of the instances are solved to optimality without preprocessing in one hour of time-limit, while 50.37\% (58.52\%) of the instances are solved to optimality using the preprocessing in one hour (10 hours) of time-limit. Furthermore, it can be seen that, even with a 10 hour time limit, many instances cannot be solved to optimality. These results are summarised in the performance profile (Figure~\ref{fig:performanceProfilesALL}) where we can observe how the formulation with preprocessing performs better than without it within one hour of time limit. Recall that the performance profile illustrates the number of solved instances over time. 

Once it has been shown that preprocessing is useful, we solve larger instances by the exact method. Specifically, we test instances with 75 nodes using the preprocessing. The results are reported in Table~\ref{tab:Exactn_75}. It is structured into three groups based on the number of facilities ($p$), the different radii ($R$), and budget values ($B_{per}$). For each group, the number of instances solved to optimality within the time limit of one hour (\#Sol.), the average solution time in seconds “t(s.)”, and the gap in percentage reported by the solver are shown.  Note that a total of 45 instances were considered for each different setting.  
\begin{table}[htbp]
  \centering
    \caption{Comparison of \namef for $n=75$, with preprocessing.}
    \medskip 
   
 \resizebox{\linewidth}{!}{
    \begin{tabular}{lrrr|lrrr|lrrr}
    $p$     & \multicolumn{1}{r}{\#Sol. } & \multicolumn{1}{r}{t(s.)} & \multicolumn{1}{r|}{\%G} & $R$     & \multicolumn{1}{r}{\#Sol. } & \multicolumn{1}{r}{t(s.)} & \multicolumn{1}{r|}{\%G} & $B_{per}$     & \multicolumn{1}{r}{\#Sol. } & \multicolumn{1}{r}{t(s.)} & \multicolumn{1}{r}{\%G} \\
    \hline
    3     &  9     & 3063.28 & 10.24 & At least one &2     & 3519.24 & 38.03& 0.025 &  13     & 2715.15 & 11.46 \\
    4     & 7     & 3042.74 & 18.25 & At least 5\% & 1     & 3572.17 & 20.38 & 0.050 &  6     & 3124.57 & 19.38 \\
    8     & 6     & 3146.23 & 33.33 & At least 10\% & 19    & 2160.85 & 3.41 & 0.100 & 3     & 3412.53 & 30.98 \\
    \hline
    \end{tabular}%
    }
   
  \label{tab:Exactn_75}%
\end{table}%
It can be seen that the number of instances solved to optimality within the time limit is small, only 22 of 135 instances, i.e., 16.3\%.  Furthermore, it can be appreciated that the majority of them correspond to the largest values of $R$. Therefore, it is interesting to develop a heuristic that yields good solutions for larger graphs. 

\subsubsection{Matheuristic solution}
The matheuristic algorithm is coded using C++ and the subproblems are solved using CPLEX 22.1.1 in Concert Technology. The experiments were conducted using an Intel(R) Xeon(R) W-2245 CPU 3.90GHz processor and 256 GB RAM. 

Intending to test the quality of the solution provided by the matheuristic,  we compute the average percentage difference between the best solution found by the heuristic algorithm ($BS_h$) and the best solution found by the exact solution method within the time limit ($BS_t$). This difference is computed as follows:
$$\%G_{BS}=\dfrac{BS_t-BS_h}{BS_t}\cdot 100.$$
If $\%G_{BS}$ is negative, then the solution found by the heuristic is better than the solution found by the exact approach. It is worth mentioning that, due to the difficulty of solving \namef, we can only provide the best solution by an exact solution method for small/medium size instances, i.e., $n=50$ and $n=75$. Therefore, $\%G_{BS}$ will be only reported for those instances.

In Table~\ref{tab:n_50}, a comparison of the different variants of the matheuristic is reported for the instances with $n=50$. This table provides the average results of 135 instances for different values of $p,$ $R,$ $B_{per},$ and the random seed, where the alternating location-downgrading search has been carried out for all eight promising sets, and afterwards one of the 1-1 local searches has been performed. The variant followed for this local search is indicated in the first row. Note that for the ``Optimal Out-In'' variant, as the procedure requires more time than other variants, we tested performing this local search once or up to ten times. For each variant, it is reported: the average solution time in seconds (t(s.)) and the GAP with the best solution found by the exact approach in one and 10 hours ($\%G_{BS}$ 1h, respectively $\%G_{BS}$ 10h). 

\begin{table}[htbp]
  \centering
      \caption{Comparison of the different variants of the matheuristic for $n=50$.}
      \medskip 
      
    \begin{tabular}{r|rr|rr|rr}
          & \multicolumn{2}{c|}{Fixed Out-In } & \multicolumn{2}{c|}{Fixed Out-Opt.\ In } & \multicolumn{2}{c}{Optimal Out-In } \\
          & a) 
          & b) 
          & a) 
          & b) 
          & One iter.
          & Up to 10 iter. 
          \\
    \hline
    t(s.)  &  0.46  & 0.59  & 1.11  & 1.08  & 2.36  & 3.15 \\
    $\%G_{BS}$ 1h & 0.86 & 0.76 & 0.34 & 0.38 & 0.03 & -0.08 \\
    $\%G_{BS}$ 10h  & 1.02& 	0.92	&0.50	&0.54	&0.20	&0.09  \\
    \hline
    \end{tabular}%

  \label{tab:n_50}%
\end{table}%
From the results shown in Table~\ref{tab:n_50} it can be seen that the matheuristic provides very good solutions. For example, the ``Optimal Out-In'' variant, in only 3 seconds, outperforms the exact approach with one-hour of time limit and provides nearly the same solutions as the exact approach with a 10 hour time limit. Furthermore, it is shown that, in general, the better the solution provided by a strategy is, the more time-consuming it is, as could be predicted in advance. Therefore, the user of the matheuristic should decide the variant to apply to find a trade-off between the quality of the solution and the required time to obtain it. Note that for the ``Fixed Out-In'' variant, it can be concluded that version b) performs better on average in these instances than version a), as it provides better solutions. Conversely, in these datasets, the ``Fixed Out-Opt. In'' variant shows that version a) achieves better results.

\begin{figure}[htbp]
   \centering
   \includegraphics[width=0.8\linewidth]{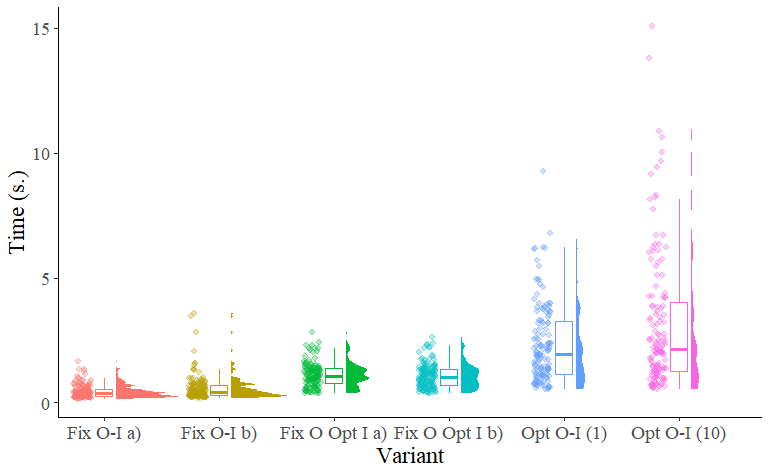}
   \caption{Comparison of time for $n=50$ by matheuristic variants.}
   \label{fig:timen_50}
\end{figure}

\begin{figure}[htbp]
   \centering
   \includegraphics[width=0.8\linewidth]{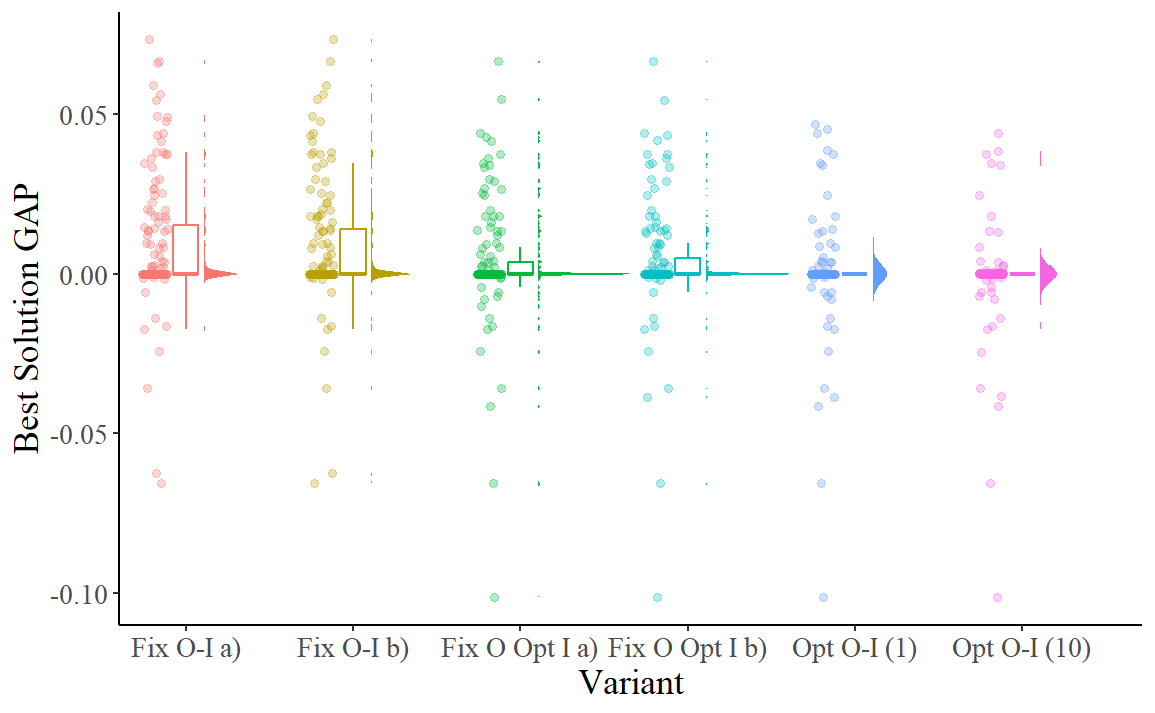}
   \caption{Comparison of $\%G_{BS}$ in one hour for $n=50$ by matheuristic variants.}
    \label{fig:GAP1n_50}
\end{figure}

To provide a comprehensive view of the data distribution and central tendency in the different variants of the matheuristic, see a rain cloud plot for the computational time in seconds (Figure~\ref{fig:timen_50}), a rain cloud plot for the $\%G_{BS}$ using one hour of time limit for the exact method (Figure~\ref{fig:GAP1n_50}). They combine elements of several different plots showing the individual data points (via the scatter plot), the statistical summaries (via the box plot), and the distribution shape (via the half-violin).

Finally, in Table~\ref{tab:n_75} a comparison of the different variants of the matheuristic is reported for the instances of size $n=75$. The structure of this table is the same as Table~\ref{tab:n_50}. The results show that the matheuristic achieves very good solutions that are better than the ones produced by the exact approach with a one-hour of time limit, and providing them for all variants in less than 20 seconds. We remark that the version up to 10 iteration of the ``Optimal Out-In'' Variant finds in less than 20 seconds 1.31\% better solutions than the exact approach in one hour. For the 135 tested instances of size $n=75$ version a) provides better solutions for variants ``Fixed Out-In'' and ``Fixed Out-Opt In''.
\begin{table}[htbp]
  \centering
      \caption{Comparison of the different variants of the matheuristic for $n=75$.}
      \medskip
    \begin{tabular}{r|rr|rr|rr}
          & \multicolumn{2}{c|}{Fixed Out-In } & \multicolumn{2}{c|}{Fixed Out-Opt.\ In } & \multicolumn{2}{c}{Optimal Out-In } \\
          & a) 
          & b) 
          & a) 
          & b) 
          & One iter.
          & Up to 10 iter. 
          \\
    \hline
     t(s.) & 1.20  & 1.26  & 2.95  & 2.77  & 8.48  & 18.94 \\
    $\%G_{BS}$ 1h & -0.16 & -0.15 & -0.61 & -0.55 & -0.90 & -1.31 \\
     \hline
    \end{tabular}%

  \label{tab:n_75}%
\end{table}%
\subsubsection{Managerial Insight}

This section tries to assess the value of considering solutions provided by our model with respect to the ones given by classical models in the literature.
In this sense, we have considered two possible attitudes of a potential location planner. In the first one, which is the optimistic one, they choose the best locations of the facilities without taking into account the existence of an attacker, i.e., they solve the MCLP in the original network $N$. Let $\ox_{UB}$ be the locations of those facilities and $S_{\ox_{UB}}$ its optimal objective value for \name, i.e., the demand that remains covered after the attacker's optimal action. The second one, on the opposite side, is the pessimistic one, where the location planner assumes that the attacker will do the maximal possible damage, i.e., they solve the MCLP in $N(u)$. Let $\ox_{LB}$ be the locations of the facilities and $S_{\ox_{LB}}$ its optimal objective value for \name. Therefore, we are comparing the bilevel model with two sequential approaches. In the sequential ones, the locations of the facilities are determined first, and then for these locations, the optimal attack is computed, i.e., \nameat is solved. Finally, the covered demand after the attack ($S_{\ox_{UB}}$ or $S_{\ox_{LB}}$) is obtained. The difference between these two attitudes with respect to our bilevel model is measured through the following two coefficients:
$$\%{MI}_O= \frac{S_{\ox_{UB}}-BS_h}{BS_{h}} \cdot 100,$$
$$\%{MI}_P= \frac{S_{\ox_{LB}}-BS_h}{BS_{h}} \cdot 100.$$
Therefore, the smaller these percentages with negative sign are, the higher the value of our model with respect to the different versions of the classical model. Indeed, they report the percentage of demand covered by the solutions of our model that is no longer covered by the ones provided by the optimistic and pessimistic attitude of the location planner, respectively.
On the other hand, these percentages also allow us to compare the quality of the matheuristic solutions obtained by the different variants for the sizes where the exact method cannot be solved.

Table~\ref{tab:MI} provides the optimistic managerial insight ($\%MI_O$) and the pessimistic managerial insight ($\%MI_P$) for the different variants of the matheuristic for the instances of size $n = 50,75,100,175,$ and $250$. The variant of the metaheuristic is indicated in the first row. Remark that the reported results are the average of 135 instances.

\begin{table}[htbp]
  \centering
    \caption{Managerial Insight of the different variants of the matheuristic.}
    \medskip 
    
    \begin{tabular}{cr|rr|rr|rr}
          $n$&& \multicolumn{2}{c|}{Fixed Out-In } & \multicolumn{2}{c|}{Fixed Out-Opt.\ In } & \multicolumn{2}{c}{Optimal Out-In } \\
          && a) 
          & b) 
          & a) 
          & b) 
          & One iter.
          & Up to 10 iter. 
          \\
    \hline
     \multirow{3}[0]{*}{50} & t(s.)  & 0.46  & 0.59  & 1.11  & 1.08  & 2.36  & 3.15  \\
   & $\%{MI}_O$  & -4.92 & -5.02 & -5.40 & -5.36 & -5.69 & -5.80  \\
    &$\%{MI}_P$ & -15.34 & -15.43 & -15.79 & -15.76 & -16.07 & -16.17  \\
    \hline
     \multirow{3}[0]{*}{75} &   t(s.) & 1.20  & 1.26  & 2.95  & 2.77  & 8.48  & 18.94 \\
   & $\%{MI}_O$ & -8.20 & -8.19 & -8.61 & -8.56 & -8.89 & -9.24 \\
    &$\%{MI}_P$ & -10.84 & -10.83 & -11.22 & -11.17 & -11.49 & -11.84  \\
    \hline
     \multirow{3}[0]{*}{100} &t(s.) & 2.31  & 2.40  & 6.13  & 5.53  & 21.03 & 58.13  \\
    &$\%{MI}_O$ & -8.74 & -8.75 & -9.01 & -8.93 & -9.23 & -9.53\\
    &$\%{MI}_P$ & -11.50 & -11.51 & -11.76 & -11.68 & -11.98 & -12.29 \\
         \hline
   \multirow{3}[0]{*}{175} &  t(s.) & 10.49 & 10.99 & 32.43 & 31.15 & 170.00 & 501.95  \\
   & $\%{MI}_O$ &-10.25 & -10.25 & -10.44 & -10.40 & -10.61 & -10.91\\
   & $\%{MI}_P$ & -5.93 & -5.93 & -6.12 & -6.08 & -6.28 & -6.59\\
     \hline
      \multirow{3}[0]{*}{250} &  t(s.) & 26.21 & 26.87 & 83.12 & 81.58 & 579.48 & 2148.67  \\
   & $\%{MI}_O$ &-14.32 & -14.32 & -14.52 & -14.46 & -14.58 & -14.92\\
    &$\%{MI}_P$ & -2.51 & -2.51 & -2.72 & -2.66 & -2.78 & -3.12 \\
    \hline
    \end{tabular}%
  
  \label{tab:MI}%
\end{table}%
From the results, it can be derived that the model and its heuristic solution are valuable, as they improve the demand covered after the attack with respect to the two considered approaches based on the classical model. The biggest improvements have achieved more than $16\%$ of additional demand coverage by the solution of our model which is no longer covered using the pessimistic approach for $n=50$ and more than $14\%$ for $n=250$ in the optimistic case.

Regarding the comparison of the different versions of matheuristics for these larger instances, the results follow the same trend in all the instance sizes, i.e., as the computational time of the matheuristic increases, the quality of the solution improves. 
Moreover, in just a few seconds, very good solutions are achieved. For these instances, there is no obvious superiority between version a) or b) for ``Fixed Out-In" variant and neither ``Fixed Out-Opt" variant. 

These managerial insight coefficients also demonstrate that different variants yield different solutions, with the ``Optimal Out-In'' providing the best, although this requires more computing time. For example, for the instances of size $n=100$, it is observed that 9.53\% (12.29\%) of the demand covered by the solution provided by our matheuristic is no longer covered by the solution of the optimistic (pessimistic) approach. Note that these matheuristic solutions are found in less than a minute on average. However, when the size of the instance is increased, it becomes much more time-consuming. 

Next, the influence of the budget on the managerial insight is examined in Table~\ref{tab:CompMangB}. In this table, the results disaggregated by $B_{per}$ for variant ``Optimal Out-In'' up to 10 iterations of the matheuristic are given. Recall that each line reports the average results of 45 instances for different values of $p$, $R,$ and random seed.

\begin{table}[htbp]
  \centering
    \caption{Effect of $B_{per}$ on the Managerial Insight for variant Optimal Out-In  Up to 10 iter.}
    \medskip
    \begin{tabular}{c|l|rrr}
    \multicolumn{1}{r}{} &       & \multicolumn{3}{c}{Optimal Out-In Up to 10 iter.} \\
\cline{2-5}    \multicolumn{1}{r}{} & $B_{per}$     & t(s.) & $\%{MI}_O$ & $\%{MI}_P$ \\
\hline 
    \multirow{3}[1]{*}{50} & 0.025 & 2.48  & -4.30 & -18.79 \\
          & 0.050 & 3.46  & -5.22 & -16.49 \\
          & 0.100 & 3.52  & -7.87 & -13.23 \\
    \hline
    \multirow{3}[2]{*}{75} & 0.025 & 16.72 & -5.97 & -14.17 \\
          & 0.050 & 20.12 & -9.09 & -12.10 \\
          & 0.100 & 19.98 & -12.66 & -9.24 \\
    \hline
    \multirow{3}[2]{*}{100} & 0.025 & 57.47 & -5.86 & -14.13 \\
          & 0.050 & 69.44 & -9.14 & -12.50 \\
          & 0.100 & 47.49 & -13.60 & -10.25 \\
    \hline
    \multirow{3}[2]{*}{175} & 0.025 & 463.80 & -6.96 & -7.97 \\
          & 0.050 & 459.80 & -10.60 & -6.74 \\
          & 0.100 & 582.26 & -15.16 & -5.04 \\
    \hline
    \multirow{3}[2]{*}{250} & 0.025 & 1409.22 & -10.17 & -4.61 \\
          & 0.050 & 1999.07 & -14.80 & -2.77 \\
          & 0.100 & 3037.72 & -19.78 & -1.96 \\
    \hline
    \end{tabular}%

  \label{tab:CompMangB}%
\end{table}%
This table shows that as the size of the dataset increases, the computational time required by the matheuristic increases. There is also a tendency for the algorithm to require less computational time as the budget gets smaller. Furthermore, it can be observed that larger values of $B_{per}$, which imply more budget allocated to downgrading the network, result in solutions whose  $\%{MI}_O$ is better than instances with less budget for downgrading. This means that the benefit of the bilevel model, with respect to the solution that does not consider an attack (optimistic approach), improves as the downgrading budget increases (observe that the optimistic approach solution will coincide with the solution of our model in case $B=0$ or $u_e=0$ for any $e\in E$). Conversely, a similar but opposite trend is observed with the managerial insight pessimistic ($\%{MI}_P$). Observe that the solution derived from the pessimistic approach will coincide with the solution of our model when $B\geq \sum_{e\in E}c_eu_e$. Although it is generally observed that the average $\%{MI}_O$ decreases and the average $\%{MI}_P$ increases with increasing $B_{per},$ there are exceptions. For example, in the instance graph50\_1 with $p=3$ and a covering radius such that each node can cover at least one, we obtain $\%{MI}_O= -6.18$, $-6.59$, and $-4.89$ as well as $\%{MI}_P=-14.93$, $-18.27$, and $-9.16$ for $B_{per}=0.025$, $0.050$, and $0.100$, respectively. 

Specifically, the benefit of considering the bilevel problem, compared to the solution of the MCLP in an auxiliary network where all edges have been maximally downgraded, is more pronounced in instances with a smaller downgrading budget. 

\begin{figure}[htbp]
    \centering
    \includegraphics[width=0.95\linewidth]{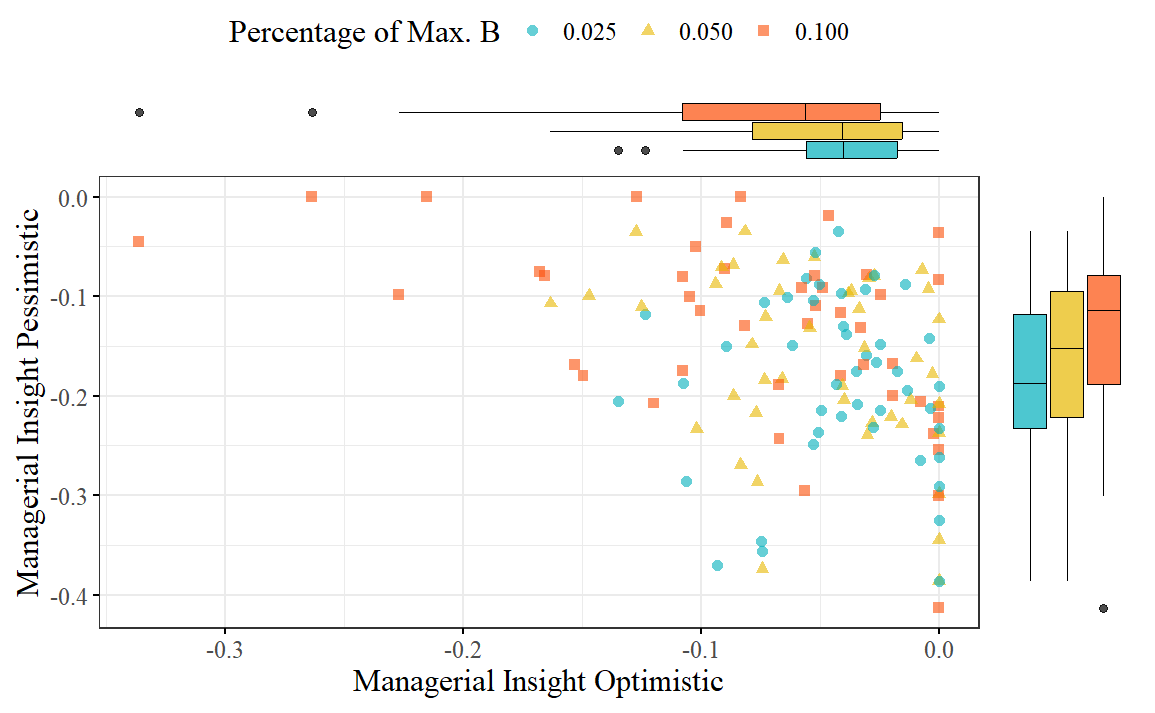}
    \caption{Comparison of Managerial Insights for $n=50$ of Optimal Out-In Variant.}
    \label{fig:MIn_50}
\end{figure}

Finally, as an example, we provide a chart (Figure~\ref{fig:MIn_50}) illustrating the effect of the downgrading budget on the managerial insight metrics for the instances of size $n=50$. For this purpose, the results of the ``Optimal Out-In'' variant, up to 10 iterations, are used.  In this graphical representation, the percentages of managerial insight optimistic $\%{MI}_O$ and pessimistic $\%{MI}_P$ are reported for each instance. Furthermore, a box plot for each value of the percentage of the maximum budget ($B_{per}$) is also presented. As before, it can be appreciated that the larger the budget, the smaller is $\% MI_O$ and the larger is $\% MI_P$, since the difference between the optimistic solution is larger and the pessimistic solution is smaller. But always, the solution of \name is better than these two approaches, for all values of $B_{per}$.  
Therefore, the main conclusions of the computational experiments are:
\begin{itemize}
\item The preprocessing improves the performance of our formulation. However, the exact approach becomes very challenging when the size of the instance is increased. 
    \item The different variants of the proposed matheuristics report quite good results in a reasonable time obtaining in a few seconds results practically identical or even better than the exact approach in one hour.
\item There is no clear superiority between the different variants of the matheuristic, as the ones that provide better solutions also require more computational time. Therefore, it should be up to the user to decide which version is best suited to their needs.
\item The value of \name has been shown, as the solution provided increases the covered demand with respect to the two straightforward sequential approaches (the solution of the MCLP in the original network and the solution of the MCLP in an auxiliary network where all the edges have been downgraded to its maximum). 

\end{itemize}

\section{Conclusion}\label{sec:Conclusion}

The downgrading maximal covering location problem is an interesting problem that addresses an open question: which is the best location for the facilities to cover the maximum demand if our network is optimally downgraded? This problem presents a paradigm in which the distance from the clients to the facilities depends on the decision variables of an external agent.

To the best of our knowledge, it is the first time this model has been presented in the literature, which has proven to be quite valuable as has been quantified in computational experiments using the managerial insight metrics. Beyond presenting a bilevel formulation and several procedures to improve it, a matheuristic has been developed that provides good solutions in a short time. Several variants for this heuristic have been proposed that help the user find the required trade-off between computational time and precision.  

This work is a starting point for several possible future research directions, such as the development of an exact bilevel method focused on this type of problem, or the consideration of other types of attacks on the network by the agent: such as considering i) a discrete number of levels of downgrading instead of being continuous,  ii) a node downgrading where the length of all incident edges is increased when a node is downgraded. 

\section*{Acknowledgements}

Research partially supported by: the Spanish Ministry of Science and Innovation through project RED2022-134149-T; Agencia Estatal de Investigación, Spain and ERDF through project PID2020-114594GB-C22; MCIN/AEI/10.13039/501100011033 and the European Union “NextGenerationEU”/PRTR through project TED2021-130875B-I00 (Marta Baldomero-Naranjo and Antonio M. Rodr\'iguez-Ch\'ia). The authors would like to thank the anonymous reviewers for their comments and
suggestions.




\appendix
\section{Appendix}\label{sec:App}
\textbf{Proof of Lemma~\ref{lemma:NPhard:AttackerProblem}}
We prove the result by reducing KNAPSACK to the decision problem \nameat-D. For a given weight threshold $T \in \mathbb{R}^+$, we define the decision version \nameat-D of \nameat as:
\begin{description}
\item[Input:] Network $N=(V,E,\ell)$, number of facilities $p$, coverage radius $R$, downgrading bounds and costs $(u_e)_{e \in E}$ and $(c_e)_{e \in E}$, respectively, budget $B$, and a set of located facilities $\overline X$.
\item[Question:] Does there exist edge downgrades $\gamma,$ such that $ 0\leq \gamma_e\leq u_e,$ for $e \in E$ with  $\sum_{e \in E} c_e \gamma_e \le B$ and $\sum_{i \in \overline C(\overline X,\gamma)}\, w_i \ge T$.
\end{description}

Let an instance of KNAPSACK be given with $n$ items of positive weight $g_i \in \mathbb{N}$ and positive value $b_i \in \mathbb{N}$, $1 \le i \le n$, knapsack capacity $K \in \mathbb{N}$, and target value $U > 0$. A solution (i.e., a Yes-Input) to KNAPSACK is a set $M \subseteq \{1,\ldots,n\}$ such that $\sum_{i \in M} g_i \le K$ and $\sum_{i \in M} b_i \ge U$. Without loss of generality, we assume $g_i \le K$, $1 \le i \le n$.

First, we observe that \nameat is in NP since a given solution for \nameat-D can be verified as such in polynomial time.
Next, given an instance for KNAPSACK, we construct an instance of the edge downgrading problem in polynomial time as follows. Let $N=(V,E,\ell)$ be a star network with central node $v_0$, \wloge we assume that the facility is located in this node, i.e., $\overline X=\{v_0\}$, 
satellite nodes $v_1,\ldots,v_n$, and edges $e_i = (v_0,v_i)$. Each satellite node $v_i$ corresponds to one item $i$, and the weight $w_i$ of the node equals the item value $b_i$. The central node is given weight $w_0=W > \sum_{i=1}^n b_i$. Finally, the weight threshold $T$ for \nameat-D is defined as $T=U$.

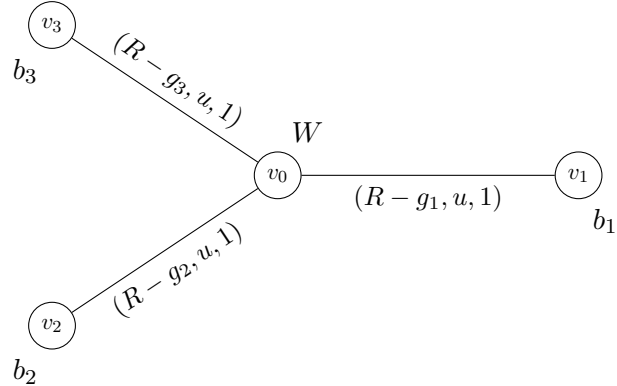
\begin{figure}[htb]
  \centering
  \begin{tikzpicture}
 	\node[scale=.8,auto=left,circle,draw,label=  80:$W$]   (v0) at ( 0, 0) {$v_0$};
 	\node[scale=.8,auto=left,circle,draw,label= -80:$b_1$] (v1) at ( 4, 0) {$v_1$};
 	\node[scale=.8,auto=left,circle,draw,label=-100:$b_2$] (v2) at (-3,-2) {$v_2$};
 	\node[scale=.8,auto=left,circle,draw,label=-100:$b_3$] (v3) at (-3, 2) {$v_3$};
 		
 	\draw (v0) -- (v1) node[pos=0.5, below, sloped]{\small $(R-g_1,u,1)$};
 	\draw (v0) -- (v2) node[pos=0.5, below, sloped]{\small $(R-g_2,u,1)$};
 	\draw (v0) -- (v3) node[pos=0.5, above, sloped]{\small $(R-g_3,u,1)$};
  \end{tikzpicture}
  \caption{Illustration for the proof of Theorem~\ref{lemma:NPhard:AttackerProblem}, with edge labels $(\ell_j,u,c)$.}
  \label{fig:complexity:star:nonuni}
\end{figure}

Let $c$ and $u$ be the unit downgrading costs and, respectively, the upper bounds on the downgrade in each edge.  We set $c=1$, $u = \max_{i=1,\ldots,n} g_i \le K,$ $R=u+1$, and the length of each edge $e_i$ as $\ell_{i} = R-g_i$. Finally, the budget equals the knapsack capacity, i.e., $B=K$. See Figure~\ref{fig:complexity:star:nonuni} for an illustration.
 
Let $M = \overline C(\overline X,\gamma)$ be the set of all satellite nodes un-covered after the downgrade $\gamma$. A satellite node $v_i$ is hereby un-covered if and only if $g_i\leq \gamma_i \leq u$. As 
\[ \sum_{i \in M}  g_i \le \sum_{i \in M}  \gamma_i \le \sum_{i=1}^n  \gamma_i \le B = K
\]
we get $\sum_{i \in M} g_i \le K$. Moreover, as 
\[ \sum_{i \in \overline C(\overline X,\gamma)} w_i =  \sum_{i \in M} w_i \ge T =  U,
\]
we obtain $\sum_{i \in M} w_i = \sum_{i \in M} b_i \ge U$ and $M$ is a solution to KNAPSACK.
Vice versa, for the same star network any solution $M$ to KNAPSACK can easily be converted into a solution for \nameat-D by setting $\gamma_i = g_i$, $i \in M$. 

As a result, the weighted \nameat-D is NP-complete on star networks and the weighted \nameat is NP-hard.

\EndProofNoNL

\bibliographystyle{abbrvnat}
\bibliography{references}

\end{document}